\documentclass[11pt]{article}

\usepackage{geometry,setspace}
\usepackage{graphicx}
\graphicspath{{Figures/}}
\usepackage{bm}
\usepackage{amssymb,amsmath,amsthm}
\usepackage{booktabs}
\usepackage{bbm}
\usepackage{caption}
\usepackage{multirow}
\usepackage{natbib}
\usepackage{authblk}

 \usepackage[normalem]{ulem}
 \usepackage[show]{ed}
 \usepackage{hyperref}
 \usepackage{cleveref}

\renewcommand{\P}{\mathbb{P}}
\newcommand{\E}{\mathbb{E}}
\newcommand{\R}{\mathbb{R}}

\newcommand{\eqdis}{\stackrel{\text{\tiny d}}{=}}

\newcommand{\abs}[1]{|#1|}

\newcommand{\var}{\operatorname{var}}

\newcommand{\cov}{\operatorname{cov}}

 \newcommand{\sign}{\operatorname{sign}}
 
\renewcommand{\leq}{\leqslant}
\renewcommand{\geq}{\geqslant}
\renewcommand{\le}{\leqslant}
\renewcommand{\ge}{\geqslant}
\newcommand{\indicator}[1]{\ensuremath{\mathbbm{1}_{\{#1\}}}}
\newcommand{\Ind}[1]{\ensuremath{\mathbbm{1}_{\{#1\}}}}
\newcommand{\powerset}{\mathcal{P}(\mathcal{D})}
\newcommand{\evenpowerset}{\mathcal{E}(\mathcal{D})}

\newtheorem{proposition}{Proposition}
\newtheorem{theorem}{Theorem}
\newtheorem{lemma}{Lemma}
\newtheorem{corollary}{Corollary}
\theoremstyle{definition}
\newtheorem{definition}{Definition}
\newtheorem{remark}{Remark}
\newtheorem{example}{Example}

\newcommand{\new}[1]{\textcolor{black}{#1}}

\begin{document}

\title{On attainability of Kendall’s tau matrices \\ and concordance signatures}

\author{A.J. MCNEIL}
\affil{The University of York Management School, University of York, UK}

\author{J.G. NE\v{S}LEHOV\'A}
\affil{Department of Mathematics and Statistics, McGill University,
  Montr\'eal, Canada}

\author{A.D. SMITH}
\affil{University College Dublin, Ireland}

\date{26 April 2022}

\maketitle

\begin{abstract}
  Methods are developed for checking and completing systems of
  bivariate and multivariate Kendall's tau concordance measures in
  applications where only partial information about dependencies
  between variables is available. The concept of a concordance signature of a multivariate
continuous distribution is introduced; this is the
  vector of concordance probabilities for margins of all orders.
It is shown that every attainable concordance signature is
  equal to the  concordance signature of a unique mixture of
  the extremal copulas, that is the copulas with extremal correlation matrices consisting
  exclusively of $1$'s and $-1$'s. A method of estimating an attainable concordance signature from data
  is derived and shown to correspond to using standard estimates of
 Kendall's tau in the absence of
  ties. The set of attainable Kendall rank correlation matrices of
  multivariate continuous distributions is proved to be identical to the set of convex combinations of extremal correlation
matrices, a set known as the cut polytope. A
methodology for testing the attainability of
  concordance signatures using linear optimization and convex
  analysis is provided.
  The elliptical copulas are shown to yield a strict
 subset of the attainable concordance signatures as well as
 a strict subset of the attainable Kendall rank correlation matrices;
 the Student $t$ copula is seen to converge, as the degrees of freedom
 tend to zero, to a
  mixture of extremal copulas sharing its concordance signature with
  all elliptical distributions that have the same correlation
  matrix. \new{A characterization of the attainable signatures of
  equiconcordant copulas is given.}
\end{abstract}


  \section{Introduction}

  In many real-world applications of statistics a modeler
  is required to impute missing information on the dependencies
  between variables, typically in the form of correlations.
This problem is particularly common in
  finance and insurance where data on certain risks
  are often sparse or non-existent. Some financial institutions use copulas
  parameterized in part by
  expert-elicited correlations to build joint models of key risks
  affecting their solvency and
  profitability~\citep{bib:embrechts-mcneil-straumann-01,bib:shaw-smith-pivak-11}.
  For example, a large insurer
  analyzing excess risk due to the Covid-19 pandemic might consider
  the interplay between mortality risk, business interruption
  risk, financial investment risk and trade credit insurance
  risk. While dependencies between these risks may be significant and
  non-negligible, they are also difficult to quantify.
  In the related area of asset management a model for the
  dependencies between asset returns is essential for optimizing a
  portfolio. However
  many assets have no track record and plausible values must
  be entered if they are to be included in the analysis. 
  Assessing dependence in the absence of data is also relevant in causal inference when unmeasured confounders are present~\citep{Stokes/Steele/Shrier:2020}.

Imputing missing information on dependence is a challenging problem because of the complex relationships between 
different pairs or subgroups of variables. The general problem of determining the
compatibility of lower-dimensional margins of higher-dimensional
distributions has only been partially resolved~\citep[see][among others]{bib:joe-96,bib:joe-97}. In this paper, we investigate the related problem of
compatibility of correlation measures for subgroups
of variables. Even for the classical linear correlation of Pearson,
the set of attainable correlation matrices has not been fully
characterized.
Because Pearson's correlation depends on the marginal distributions,
the set of attainable matrices is generally much more complicated than
the elliptope of positive semi-definite correlation matrices
investigated, for example, by~\citet{bib:huber-maric-15,bib:huber-maric-19},
unless attention is restricted to special choices of margins such as the normal \citep{bib:embrechts-mcneil-straumann-01}.   Recently, \citet{bib:hofert-koike-19} described the set of attainable Blomqvist's
beta matrices, while
\citet{Embrechts/Hofert/Wang:2016} found conditions characterizing 
matrices of tail dependence coefficients. \citet{Devroye/Letac:2015}
and~\citet{Wang/Wang/Wang:2019} consider the set of attainable
Spearman rank correlation matrices; this coincides with the elliptope
of linear $d\times d$ correlation matrices when $d \le 9$ but not when $d \ge 12$,
while the case where $d\in\{10,11\}$ remains to be settled.

In this paper, we provide a complete solution to the
attainability and compatibility problem when dependence is measured by the widely-used
Kendall's tau rank correlation~\citep{Kendall:1938,bib:kruskal-58,bib:joe-90} which parametrizes many popular dependence models. In order
to take into account bivariate and higher-order associations
between subsets of variables, we quantify the dependence of a
multivariate random vector using a vector-valued measure that we call
the concordance signature, which underlies all bivariate and
multivariate Kendall's tau coefficients; this concept is made precise
in Section~2.

\new{A note on our use of the terms attainability and compatibility may be helpful at this point. We use attainable
when we talk about a logically coherent collective entity that could belong to
a probability distribution, such as a
concordance signature or a Kendall's tau matrix. We use compatible when we talk about the relationship
between sub-components of an attainable entity, such as Kendall's tau rank correlations for
different pairs of variables. We also refer to attainable signatures as
being compatible with probability distributions.}

Our main result is to fully characterize the set of attainable concordance signatures of continuous multivariate distributions. As a by-product, we prove the
conjecture of~\citet{bib:hofert-koike-19} that the set of attainable
Kendall rank correlation matrices is identical to the set of convex combinations of the extremal correlation
matrices, i.e., the correlation matrices consisting exclusively of
$1$'s and $-1$'s; this set is also known as the cut
polytope~\citep{bib:laurent-poljak-95}.  As we show, this characterization follows from the links between concordance signatures, multivariate
Bernoulli distributions, and extremal mixture copulas, i.e., mixtures of the $2^{d-1}$ possible copulas with extremal correlation
matrices~\citep{bib:tiit-96}. We explain why the set of attainable
Kendall correlation matrices is identical to the set of attainable correlation
matrices for multivariate Bernoulli random vectors with symmetric margins as derived
by~\citet{bib:huber-maric-15,bib:huber-maric-19}.

\begin{figure}[t!]
  \centering
  \includegraphics[width=9cm,height=6.5cm]{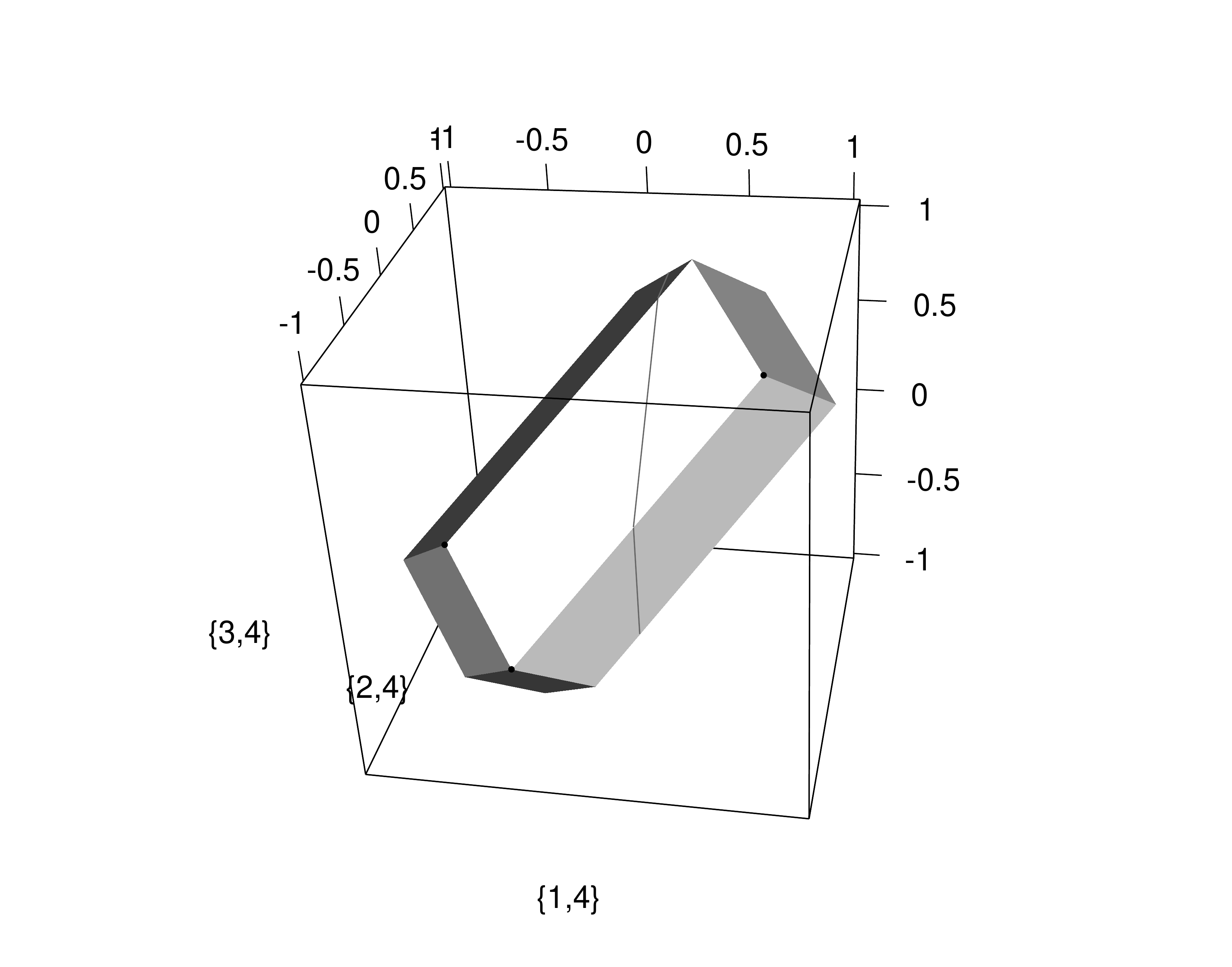}
  \includegraphics[width=5.5cm,height=6.5cm]{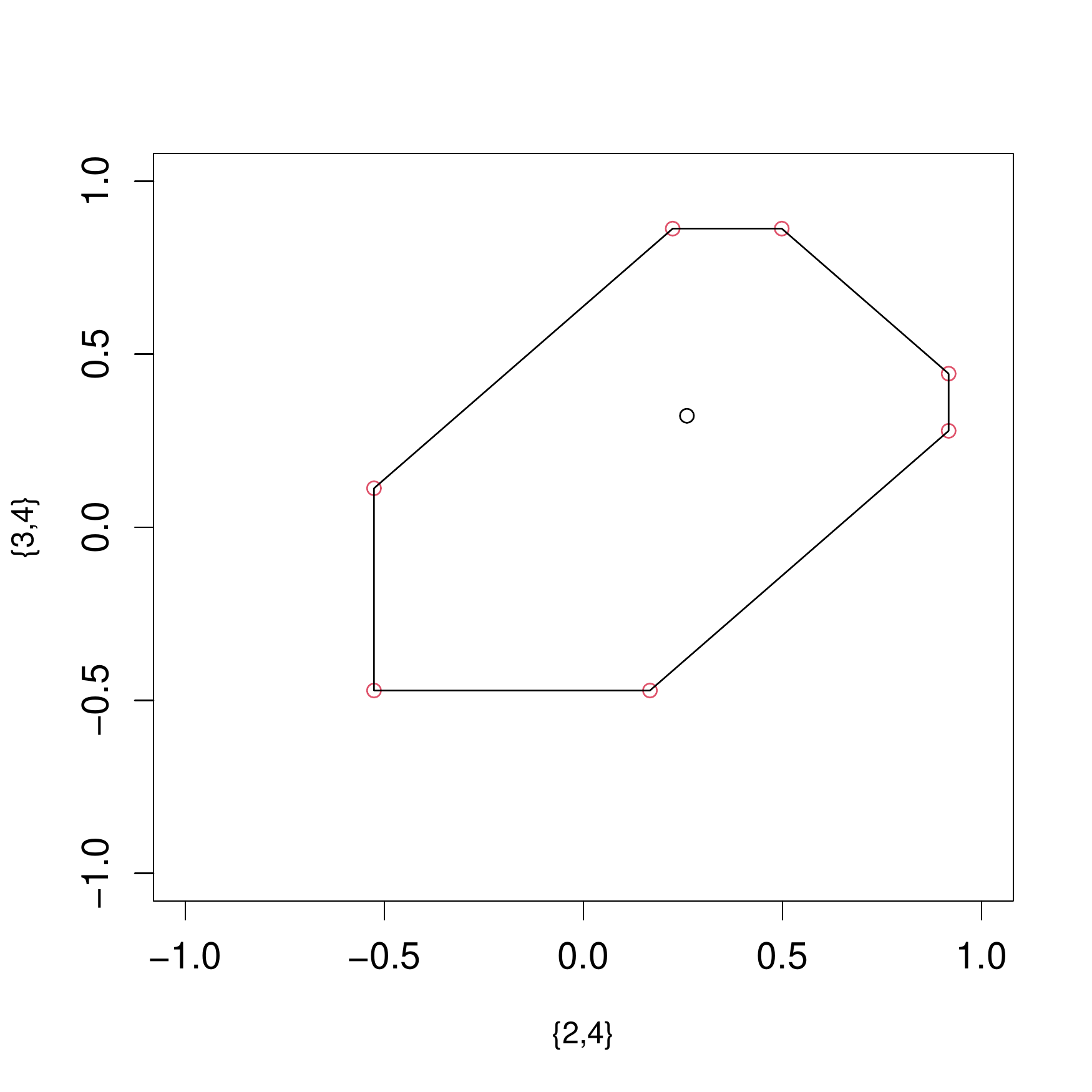}
 \caption{Left panel: set of attainable Kendall rank correlations
   for X4coin paired with three other cryptocurrencies in
   Example~\ref{example-data}. Right panel: set of attainable
   Kendall's tau values for (X4coin, Etherium) and (X4coin, Litecoin)
   when the (X4coin, Bitcoin) value is fixed at 0.598.}
  \label{fig:estimation}
\end{figure}

The methodology we propose has a number of
important applications. First, it allows us to determine whether a set of estimated or
expert-elicited Kendall correlations (bivariate or multivariate) is compatible with a valid multivariate distribution. If it is, we then propose a method of determining the set in which any remaining unmeasured Kendall correlations must lie using standard techniques from linear optimization and convex analysis. To illustrate, suppose, for example, we have Bitcoin, Etherium and Litecoin in our
  portfolio and  we invest in a new cryptocurrency `X4coin'.
  Based on return data for Bitcoin, Etherium and Litecoin in 2017, the attainable Kendall rank correlations between the as
  yet unobserved X4coin and
  the other three currencies can be shown to lie in the
  three-dimensional polytope in the left panel of
  Figure~\ref{fig:estimation}. As soon as we form an opinion on the
  rank correlation between X4coin and one of the cryptocurrencies,
  this further restricts the set of possible rank correlations
  between X4coin and the other two currencies, as shown in the right panel of
  Figure~\ref{fig:estimation}; full details are given later in Example~\ref{example-data}.
  
We also prove that sample concordance signatures based on classical estimators of 
Kendall's tau are themselves concordance signatures of valid multivariate distributions. Our finding that the concordance signature of a
multivariate distribution with continuous margins is always equal to the concordance
signature of a unique mixture of the extremal copulas then offers a powerful technique in Monte
Carlo simulation or risk analyses. It allows us to draw realizations from a model with identical concordance signature
to that estimated from the data
but with an extreme form of tail~dependence.

Finally, our methodology allows us to take a closer look at the dependence structure
inherent in different classes of distributions.  We show that the concordance signatures of the family of
elliptical copulas form a strict subset of the set of  all possible attainable concordance
signatures. The surprising consequence is the existence
of Kendall correlation matrices that do not
correspond to any elliptically distributed random vector; this behavior
provides a contrast with Pearson correlation matrices, each of which is the correlation matrix of at least one elliptical distribution. Moreover, we show that the 
$d$-dimensional Student $t$ copula with
correlation matrix $P\in
\R^{d\times d}$ and degree-of-freedom parameter $\nu>0$ converges
pointwise to a mixture of extremal copulas which shares its
concordance signature as $\nu\to 0$. 
\new{We also investigate the concordance signatures of exchangeable and
equiconconcordant distributions.}


\section{Concordance signatures} \label{sec:2}

Throughout this paper, we use bold symbols such as
$\bm{x}=(x_1,\dots,x_d)$ to denote vectors in $\mathbb{R}^d$ and
understand expressions such as $\bm{x}+\bm{y}$ as componentwise
operations; similarly $\bm{x} \leq \bm{y}$
implies that all components are ordered.

 Let $\bm{X}=(X_1,\ldots,X_d)$ denote a generic random vector with continuous
 margins $F_1,\ldots,F_d$. As is well known, the unique
 copula $C$ of $\bm{X}$ is the distribution function of the  random vector
 $\bm{U}=(U_1,\ldots,U_d)$ where $U_i = F_i(X_i)$ for
 $i=1,\ldots,d$ \citep{bib:joe-97,Nelsen:1999}.

 For any subset $I
 \subseteq \mathcal{D}=\{1,\ldots,d\}$ with $I \neq \emptyset$ we write $\bm{X}_I$ and $\bm{U}_I$ to denote
 sub-vectors of $\bm{X}$ and $\bm{U}$, and $C_I$ for the copula of
 $\bm{X}_I$, i.e., the distribution function of $\bm{U}_I$. The multivariate probability of concordance of $\bm{X}_I$ is defined as
 \begin{equation}
   \label{eq:1}
   \kappa_I = 
     \P\left(\{\bm{X}_I \leq
     \bm{X}^*_I\} \cup \{\bm{X}^*_I \leq
     \bm{X}_I\} \right),
 \end{equation}
where $\bm{X}^*$ is a random vector independent of $\bm{X}$ but with the same distribution. Note that~\eqref{eq:1} implies $\kappa_{\{i\}} = 1$ for a singleton
and we adopt the convention $\kappa_\emptyset = 1$. 

For subsets with cardinality $|I| \geq 2$, the concordance
probabilities quantify the association between the components of
$\bm{X}_I$. Indeed $\kappa_I$ is
related to Kendall's tau $\tau_I$ \citep{Kendall:1938, bib:kruskal-58} when $|I|=2$  and to its multivariate analogue  \citep{bib:joe-90,bib:genest-neslehova-ben-ghorbal-11} when $|I|>2$ via the formula
\begin{equation}
  \label{eq:14}
(2^{|I|-1}-1) \tau_I = 2^{|I|-1}\kappa_I -1 \;.
\end{equation}
Moreover, following \citet{Nelsen:1999}, $\kappa_I$ depends only on the copula $C_I$ of $\bm{X}_I$
according to
\begin{equation}\label{eq:13}
  \kappa_I  = 2\P\left(\bm{X}_I \leq
    \bm{X}^*_I\right)=2\P\left(\bm{U}_I \leq
    \bm{U}^*_I\right) =2\int_{[0,1]^{|I|}}C_I(\bm{u}) d C_I(\bm{u}).
\end{equation}
Without loss of generality, we can thus refer to $\kappa_I$ as the
concordance probability of the copula $C_I$ of $\bm{X}_I$ and use the
notation $\kappa_I = \kappa(C_I)= \kappa(\bm{X}_I)= \kappa(\bm{U}_I)$ interchangeably.

To summarize the various dependencies inherent in the vector $\bm{X}$, we can
collect the values $\kappa_I$ for all subsets $I \in \mathcal{P}(\mathcal{D})$,
where $\mathcal{P}(\mathcal{D})$ denotes the power set of
$\mathcal{D}$. This motivates calling the vector $( \kappa_I: I \in
\mathcal{P}(\mathcal{D}))$ the {\it full concordance signature} of
$\bm{X}$. Sub-vectors of the full concordance signature are known as
{\it partial concordance signatures}. From Proposition~1 in \citet{bib:genest-neslehova-ben-ghorbal-11}, which is based on the exclusion-inclusion principle, we can deduce that for any set $I\subseteq\mathcal{D}$ of odd cardinality,
\begin{equation}\label{eq:oddkappa}
\kappa_I  = 1 + \sum_{A \subset I, 1 \le |A| < |I|} (-1)^{|A|} \P(\bm{U}_A \leq \bm{U}^*_A) = 1 - \frac{|I|}{2} +  \sum_{A \subset I, 2 \le |A| < |I|} (-1)^{|A|}  \frac{\kappa_A}{2} \; .
\end{equation}
Thus the full concordance signature can be deduced from the concordance probabilities for
subsets of $\mathcal{D}$ of even cardinality. We formalize ideas in the following definition.
%
%

\begin{definition}\label{def:1}
 A label set is any collection $S$ of subsets of $\mathcal{D}$ such that
    $\emptyset \in S$.  The vector $\bm{\kappa}_S(C) =(\kappa(C_I)  :
    I \in S)$ is called the partial concordance signature of the copula
    $C$ for the label set $S$. By
     convention, the elements of $S$ are taken in
     lexicographical order.
 When the label set is the even power set $\evenpowerset = \{I: I
\subseteq \mathcal{D}, |I|\,\text{even}\}$ (containing the empty set $\emptyset$ by convention) 
 the partial concordance
         signature is called the even
       concordance signature of $C$ and is denoted $\bm{\kappa}(C)$.
 When the label set is $\powerset$ the partial concordance
         signature is called the full
       concordance signature of $C$ and is denoted $\tilde{\bm{\kappa}}(C)$.
\end{definition}

\section{Extremal mixture copulas}\label{sec:notat-gener-case}

Our main tools for the study of concordance signatures are mixtures of so-called extremal copulas and their relationship with multivariate
Bernoulli distributions. In this section we develop the necessary
notation and theory.

An extremal copula $C$ with index set $J \subseteq \mathcal{D}$ is the distribution function of the random vector $\bm{U}=(U_1,\ldots, U_d)$ where for $j \in \mathcal{D}$, $U_j \eqdis  U$ if $j \in J$, $U_j \eqdis 1-U$ if $j \in J^\complement$,
and $U$ is a standard uniform random variable. For all
$\bm{u} \in [0,1]^d$ it has the explicit form
\begin{equation}\label{eq:4}
  C(\bm{u}) = ( \min_{j \in J} u_j + \min_{j\in J^\complement}
    u_j  -1)^+\;,
\end{equation}
where for any $x \in \mathbb{R}$, $x^+ = \max(x,0)$ denotes the
positive part of $x$ and we employ the convention  $\min_{j \in \emptyset} u_j
= 1$.

In dimension $d\geq 2$ there are $2^{d-1}$ extremal copulas and we enumerate
them in the following way. For $k\in \{1,\ldots,2^{d-1}\}$ let $\bm{s}_{k}=(s_{k,1},\ldots,s_{k,d})$ be the
vector consisting of the digits of $k-1$ when represented
  as a $d$-digit binary number.
For example, when $d=4$ we have exactly eight extremal copulas
corresponding to the vectors
\begin{displaymath}
\begin{aligned}
  \bm{s}_{1} &= (0,0,0,0),\quad \bm{s}_{2} = (0,0,0,1),\quad \bm{s}_{3} =
  (0,0,1,0),\quad \bm{s}_{4} = (0,0,1,1), \\
\bm{s}_{5} &= (0,1,0,0),\quad \bm{s}_{6} =(0,1,0,1),\quad \bm{s}_{7}=(0,1,1,0),\quad \bm{s}_{8}=(0,1,1,1).
\end{aligned}
\end{displaymath}
For each $k\in \{1,\ldots,2^{d-1}\}$, $\bm{s}_k =(s_{k,1},\ldots,s_{k,d})$ and $\bm{1}-\bm{s}_k=(1-s_{k,1},\ldots,1-s_{k,d})$ are opposite vertices of the unit hypercube connected by one of its $2^{d-1}$ main diagonals, $\bm{1}$ being the vector of ones. The $k$th extremal copula in dimension $d$, denoted $C^{(k)}$, spreads its probability mass uniformly along the latter diagonal joining $\bm{s}_k$ and $\bm{1}-\bm{s}_k$.  Its index set $J_k \subseteq \mathcal{D} = \{1,\ldots,d\}$ is defined 
as the set of indices corresponding to zeros in $\bm{s}_k$; that is, $j \in J_k$
if $s_{k,j} = 0$ and $j \in J_k^\complement $ if $s_{k,j} = 1$. For the 4-dimensional
example above we have $J_1 = \{1,2,3,4\}$, $J_1^\complement = \emptyset$, $J_2 =
\{1,2,3\}$, $J_2^\complement = \{4\}$ and so on. 
Note also that $J_1 = \mathcal{D}$ so that $C^{(1)}(\bm{u}) =\min(u_1,\ldots,u_d)$ is the
 comonotonicity or Fr\'echet--Hoeffding upper bound copula.

Extremal copulas owe their name to the fact that their correlation
matrices are extremal correlation matrices, that is, matrices consisting
exclusively of 1's and $-1$'s. Indeed, for any $k \in \{1,\ldots,
2^{d-1}\}$, we find that the matrix of pairwise Kendall's tau values for the $k$th extremal copula $C^{(k)}$  is $P^{(k)} =
(2\bm{s}_k -\bm{1}) (2\bm{s}_k -\bm{1})^\top$.  This matrix is simultaneously
the matrix of pairwise Pearson and Spearman correlations of $C^{(k)}$.

   
  
  \begin{definition}
  An extremal mixture copula is a copula of the form $C^* = \sum_{k=1}^{2^{d-1}} w_{k} C^{(k)}$, where for all $k \in\{1,\ldots, 2^{d-1}\}$, $C^{(k)}$ is the $k$th extremal copula, $w_k \geq 0$, and $\sum_{k=1}^{2^{d-1}} w_k =1$.
  \end{definition}
 

The following proposition shows how the extremal mixture copulas are
related to multivariate Bernoulli distributions; the proof is given in~\ref{app:A}.
\begin{proposition}\label{prop:bernoulli}
Let $U$ be a standard uniform random variable and $\boldsymbol{B}$ a
$d$-dimensional multivariate Bernoulli vector independent of $U$. Then the distribution function of the vector 
 \begin{equation}
\label{eq:7}
U \boldsymbol{B} + (1-U)(\boldsymbol{1} -\boldsymbol{B})
 \end{equation}
 is an extremal mixture with weights given, for each $k \in \{1,\ldots, 2^{d-1}\}$, by
\begin{equation}\label{eq:8}
w_k = \P(\boldsymbol{B} = \boldsymbol{s}_k) +  \P(\boldsymbol{B} = \bm{1}-\boldsymbol{s}_k),
\end{equation}
where $\boldsymbol{s}_k$ is the vector consisting of the digits of $k-1$ when represented as a $d$-digit binary number. Conversely, any extremal mixture copula $C^* = \sum_{k=1}^{2^{d-1}} w_k C^{(k)}$  is the distribution function of a random vector of the form~\eqref{eq:7}, where $U$ is independent of $\boldsymbol{B}$ and~\eqref{eq:8} holds for all $k \in \{1,\ldots, 2^{d-1}\}$.
\end{proposition}

The class of Bernoulli distributions satisfying~\eqref{eq:8} is
infinite since the mass $w_k$ can be split between the events $\{\boldsymbol{B}=\boldsymbol{s}_k\}$ and
$\{\boldsymbol{B}=\bm{1}-\boldsymbol{s}_k\}$ in an arbitrary way. It is
important for the arguments used in this paper to single out a representative and we do this by setting
$\boldsymbol{B}=\boldsymbol{Y}$, where $\boldsymbol{Y}$ is \textit{radially
symmetric} about $\boldsymbol{0.5}$. This means that $\bm{Y}
\stackrel{d}{=} \bm{1} - \bm{Y}$
and implies that, for each $k \in
\{1,\ldots, 2^{d-1}\}$,
\begin{equation}
  \label{eq:41}
  \P(\boldsymbol{Y} = \boldsymbol{s}_k) =  \P(\boldsymbol{Y} = \boldsymbol{1}-\boldsymbol{s}_k) = 0.5 w_k.
\end{equation}
Such distributions are
also known as {\it palindromic} Bernoulli
distributions~\citep{Marchetti/Wermuth:2016} and they are fully
parameterized by the $2^{d-1}$ probabilities $\P(\bm{Y} = \bm{s}_k) = 0.5 w_k$, $k \in \{1,\ldots, 2^{d-1}\}$.
This means that there is a bijective mapping
between the extremal mixture copulas and the radially symmetric multivariate Bernoulli
distributions in dimension $d$. 
A number of constraints apply to radially symmetric Bernoulli random vectors. Apart from the fact that  $\P(Y_i =1) = 0.5$ for all $i$, we will make use of the following insight, proved in~\ref{app:A}.

\begin{proposition}\label{prop:type-2-repr}
  \begin{enumerate}
  \item[(i)] The radially symmetric Bernoulli distribution of a vector $\bm{Y}$ in
 dimension $d$ is uniquely determined by the vector of probabilities
  $\bm{p}_{\bm{Y}} = (p_I  :
  I \in \evenpowerset \setminus\emptyset)$, where $p_I = \P(\bm{Y}_I = \bm{1})$ and $\evenpowerset$ is the even power set as in Definition~\ref{def:1}.
  \item[(ii)] The vector
    $\bm{p}_{\bm{Y}}$ in (i) is the shortest vector of the form  $ (p_I  :
  I \in S\setminus\emptyset)$ for a label set $S \subset \powerset$ which uniquely determines the
  distribution of $\bm{Y}$ for all \new{radially symmetric Bernoulli} random vectors $\bm{Y}$.
  \end{enumerate}
\end{proposition}

We close this section by noting that the independence
between $U$ and $\boldsymbol{Y}$ in the stochastic representation
$\bm{U} \eqdis U\bm{Y} + (1-U)(\bm{1}-\bm{Y})$ of an extremal mixture
copula cannot be dispensed with. An example of a copula that violates this condition is provided in~\ref{app:A} and reveals an interesting contrast between the extremal
copulas and the mixtures of extremal copulas. While a necessary and
sufficient condition for a vector $\bm{U}$ to be distributed according
to an extremal copula is that its bivariate margins should be extremal
copulas \citep{bib:tiit-96}, the analogous statement does not hold for extremal mixture
copulas. Although it is necessary that the bivariate margins of
an extremal mixture copula are extremal mixture copulas, the
example shows that this is not sufficient. An additional condition is required, as detailed in the next proposition which is proved in~\ref{app:A}.

\begin{proposition}\label{theorem:pairwise}
  The distribution of a random vector $\bm{U}=(U_1,\ldots,U_d)$ is a mixture of
  extremal copulas if and only if its  bivariate marginal
  distributions are mixtures of extremal copulas and for all $u \in [0,1]$,
  \begin{equation}\label{eq:10}
    \P\left(U_1 \leq u \mid \indicator{U_j = U_1}, j\neq 1\right) =
    u\;.
  \end{equation}
\end{proposition}


\section{Characterization of Concordance Signatures}\label{sec:mult-prob-conc}

In this section we characterize concordance signatures of
arbitrary copulas in dimension~$d$. To do so, we first calculate
concordance signatures of extremal mixture copulas and
investigate their properties. For any  $k\in\{1,\ldots,2^{d-1}\}$ let
$C^{(k)}$ be an extremal copula with index set $J_k$ as specified in Section~\ref{sec:notat-gener-case}. For $I
\subseteq \mathcal{D}$ we introduce the notation
\begin{equation}\label{eq:6}
  a_{I,k} =
  \begin{cases}
    1&\text{if $I \subseteq J_k$ or $I \subseteq J_k^\complement$,}\\
    0&\text{otherwise,}
    \end{cases}
  \end{equation}
noting that $a_{\{i\},k} = 1$ and $a_{\emptyset,k}
  = 1$ for all $k$ and $i$. Then we have the following result.

\begin{proposition}\label{prop:multi-concordance}
\begin{enumerate}
\item[(i)]  For any $k \in \{1,\ldots, 2^{d-1}\}$ and $I \in \powerset$, $\kappa(C^{(k)}_I) =
a_{I,k}$.
\item[(ii)] If $C$ is an extremal mixture copula of the form  $\sum_{k=1}^{2^{d-1}} w_{k} C^{(k)}$ then, for any $I \subseteq\mathcal{D}$,
\begin{equation}\label{eq:15}
  \kappa(C_I) = \sum_{k=1}^{2^{d-1}} w_{k}
  \kappa\left(C^{(k)}_I\right) = \sum_{k=1}^{2^{d-1}} w_{k}a_{I,k}\; .
\end{equation}
\end{enumerate}
\end{proposition}
\begin{proof} 
To show (i), it suffices to consider sets with $|I| \ge 2$. 
 For
$k\in\{1,\ldots,2^{d-1}\}$, let
$\bm{X}^{(k)}=(X^{(k)}_{1},\ldots,X^{(k)}_{d})$ be a random vector
with continuous margins and copula $C^{(k)}$. Then the sets $\{X^{(k)}_{j} : 
j \in J_k\}$ and $\{X^{(k)}_{j}: 
j \in J_k^\complement\}$ are sets of comonotonic random variables which are concordant with probability $1$, while any pair
$(X^{(k)}_{i}, X^{(k)}_{j})$ such that $i\in J_k$ and $j \in J_k^\complement$, or vice versa, is
a pair of countermonotonic random variables and hence 
discordant with probability $1$.

To establish (ii), note again that~\eqref{eq:15} holds trivially for sets $I$ which are
  singletons or the empty set. For sets such that $|I|\geq 2$ we can
  use~\eqref{eq:13} to write
  \begin{displaymath}
\kappa(C) = 2  \int_{[0,1]^{d}}C(\bm{u}) d C(\bm{u}) =
    2\sum_{j=1}^{2^{d-1}} \sum_{k=1}^{2^{d-1}} w_j
    w_k \int_{[0,1]^{d}} C^{(j)}(\bm{u}) d C^{(k)}(\bm{u})\;.
  \end{displaymath}
  Introducing independent random vectors $\bm{U}^{(j)}\sim C^{(j)}$ and
  $\tilde{\bm{U}}^{(k)} \sim C^{(k)}$ for $j
 \in\{1,\ldots, 2^{d-1}\}$ and $k\in\{1,\ldots, 2^{d-1}\}$ we calculate that
  \begin{displaymath}
    \int_{[0,1]^{d}} C^{(j)}(\bm{u}) d C^{(k)}(\bm{u}) =
    \P\left( \bm{U}^{(j)} \leq \tilde{\bm{U}}^{(k)} \right) =
    \begin{cases}
      \frac{1}{2}& \text{if $j=k=1$,}\\
      \frac{1}{4}&\text{if $j=1$ or $k=1$ but $j\neq k$,}\\
      0&\text{if $j\neq 1$ and $k \neq 1$.}
      \end{cases}
  \end{displaymath}
  Hence we can verify that
$
    \kappa(C) = 2 (w_1^2/2) + 4 w_1\{ (w_2/4) +
      \cdots + (w_{2^{d-1}}/4)\} = w_1^2 + w_1(1- w_1) = w_1,
$
which
  is the weight on $C^{(1)}
  = M$, the $d$-dimensional comonotonicity copula. If $\bm{U}$ is
  distributed according to $C$
  then the vector
  $\bm{U}_I$ is distributed according to
a mixture of extremal copulas in dimension $|I|$ and
  it follows that  $\kappa(\bm{U}_I) = \kappa(C_I) = \tilde{w}_1$ where
  $\tilde{w}_1$ is the weight attached to the case where
  $\bm{U}_I$ is a comonotonic random vector. This is given by
  \begin{displaymath}
    \tilde{w}_1 = \sum_{k=1}^{2^{d-1}} w_{k}a_{I,k} = \sum_{k=1}^{2^{d-1}} w_{k} \kappa\left(C^{(k)}_I\right), 
  \end{displaymath}
  where the second equality follows from part (i).
\end{proof}




\new{
\begin{remark}
It may be noted that the extremal copulas $C^{(k)}$ are the only
copulas that give rise to extremal concordance signatures, that is
concordance signatures consisting only of zeros and ones as in part
(i) of Proposition~\ref{prop:multi-concordance}. This follows from the
fact that two variables are comonotonic if and only if the bivariate
concordance probability is one and countermonotonic if and only if the
bivariate concordance probability is zero \citep{bib:embrechts-mcneil-straumann-01}; thus a copula with an
extremal signature must have bivariate margins that are extremal
copulas. \citet{bib:tiit-96} showed that, if the bivariate margins of
  a copula are extremal copulas, then the copula is an extremal copula. 
\end{remark}
}

%

From Proposition~\ref{prop:multi-concordance} we see that for any label set $S$ and any $k \in \{1,\ldots,2^{d-1}\}$, the partial concordance signature of $C^{(k)}$ is $\bm{\kappa}_S(C^{(k)})=(a_{I,k} : I \in
S)$. We also see that the partial concordance signature of an extremal
mixture copula $C=\sum_{k=1}^d w_k C^{(k)}$ is a convex combination of
partial signatures of extremal copulas with the same weights, so that
$\bm{\kappa}_S(C) = \sum_{k=1}^d w_k \bm{\kappa}_S(C^{(k)})$. The following two key properties of the even concordance signature of an
  extremal mixture link back to the Bernoulli
  representation in Proposition~\ref{prop:bernoulli}; the proof relying on Proposition \ref{prop:type-2-repr} is given in \ref{app:A}.
  
 \begin{proposition}\label{prop:main-result-prelim}
    Let $\bm{a}_k =  (a_{I,k} : I \in \evenpowerset)$, for $k \in
    \{1,\ldots,2^{d-1}\}$, and let $\bm{\kappa}(C) =
      \sum_{k=1}^{2^{d-1}} w_k \bm{a}_k$ be the even concordance
      signature of the
      extremal mixture $C=\sum_{k=1}^d w_k C^{(k)}$.
    \begin{enumerate}
    \item[(i)] $\bm{\kappa}(C) $ 
    uniquely determines $C$. Moreover, it is the minimal partial
      concordance signature of $C$ which uniquely determines $C$ in
      all cases.
      \item[(ii)] The vectors $\bm{a}_k$, $k \in \{1,\ldots,2^{d-1}\}$ are
      linearly independent.
    \end{enumerate}

  \end{proposition}

  \begin{proof}
    For part (i) recall that a random vector $\bm{U} \sim C$ has the stochastic
    representation
    $\bm{U} \eqdis U\bm{Y} + (1-U)(\bm{1}-\bm{Y})$ where $\bm{Y}$ is a random vector
    with a radially symmetric Bernoulli distribution. For any
    set $I \in \evenpowerset\setminus\emptyset$ the
    components of $U\bm{Y} + (1-U)(\bm{1}-\bm{Y})$ are concordant if and only if $\bm{Y}_I
    =\bm{1}$ or $\bm{Y}_I = \bm{0}$. It follows from the radial
    symmetry property that $\kappa_I(C) = 2\P(\bm{Y}_I =\bm{1})$. By
    Proposition~\ref{prop:type-2-repr} the vector
    $(\P(\bm{Y}_I = \bm{1}) : I \in
    \evenpowerset\setminus\emptyset)$ is the minimal vector of event
    probabilities that pins down the law of $\bm{Y}$ in all cases and hence, by
    equation~\eqref{eq:41}, the
    weights $w_k = 2\P(\bm{Y} = \bm{s}_k)$ in the representation $C =
    \sum_{k=1}^{2^{d-1}}w_k C^{(k)}$. If another extremal
    mixture copula
    $\tilde{C}$ shares the same values for the concordance
    probabilities of even order, then the weights must be identical.

    For part (ii) suppose, on the contrary, that the vectors $\bm{a}_k$ are linearly dependent, that is, there exist scalars $\mu_k$, $k \in \{1,\ldots, 2^{d-1}\}$ such that  $\mu_k \neq 0$ for at least one $k$ and 
\begin{equation}\label{eq:B1}
\sum_{k=1}^{2^{d-1}} \mu_k \bm{a}_k=\bm{0}. 
\end{equation}
 To prove the result, we need to show that this assumption leads to a
 contradiction. We select a copula $C$ such that $C=
 \sum_{k=1}^{2^{d-1}} w_k C^{(k)}$ and $w_k >0$ for all $k$.
Because the first component of $\bm{a}_k$ is 1 for all $k$, we have that $\sum_{k=1}^{2^{d-1}} \mu_k =0$. Hence, there exists at least one $k^*$ such that $\mu_{k^*} > 0$. From \eqref{eq:B1} we have that for each $\alpha \in \mathbb{R}$, 
 $\sum_{k=1}^{2^{d-1}} (w_k - \alpha \mu_k) \bm{a}_k =\bm{\kappa}(C)$.
 
 Let $\alpha^* = \min ( (w_k/\mu_k) \: : \; \mu_k > 0)$. Clearly,
 $\alpha^* > 0$ because the set over which the minimum is taken
 contains at least $w_{k^*} / \mu_{k^*}$ and because $w_k > 0$ for all
 $k$. Now define $\bm{w}^* = \bm{w} - \alpha^* \bm{\mu}$. These
 weights are non-negative and sum up to one, while $\bm{w}^* \neq
 \bm{w}$. This implies that $\bm{\kappa}(C) = \sum_{k=1}^{2^{d-1}}
 w_k \bm{a}_k =  \sum_{k=1}^{2^{d-1}} w_k^* \bm{a}_k $ but this is not
 possible because the mixture weights are unique by part (i).
\end{proof}

Proposition~\ref{prop:main-result-prelim} shows that the set
of attainable even concordance signatures of $d$-dimensional extremal
mixture copulas is the convex hull
$$
\mathbb{K}=\left\{ \sum_{k=1}^{2^{d-1}} w_k \bm{a}_k \: : \: w_k \ge 0, k=1,\ldots, 2^{d-1}, \sum_{k=1}^{2^{d-1}} w_k =1\right\}
$$
and implies in particular that
$\mathbb{K}$ is a convex polytope with vertices $\bm{a}_k$, $k \in
\{1,\ldots, 2^{d-1}\}$. Similarly the set of full concordance
signatures of extremal mixtures is the convex hull of $\tilde{\bm{a}}_k=(a_{I,k} : I \in \powerset)$,  $k \in
\{1,\ldots, 2^{d-1}\}$. 

We are now ready for the main result of this paper which shows
that these sets are also the sets of attainable even and full concordance
signatures for any $d$-dimensional copula, and hence of any random vector with continuous margins.

\begin{theorem}\label{theorem:conc-sign-copul}
Let $C$ be a $d$-dimensional copula and $\bm{\kappa}(C) = (\kappa_I:
I \in \evenpowerset )$ its even concordance signature. Then
there exists a unique extremal mixture copula $C^* =
\sum_{k=1}^{2^{d-1}}w_k C^{(k)}$ such that
$\bm{\kappa}(C) = \bm{\kappa}(C^*)$ . The weights $w_k$, which are non-negative and add up to $1$, are the unique solution to the system of $2^{d-1}$ linear equations given by 
\begin{equation}\label{eq:34}
    \bm{\kappa}(C)
   = \sum_{k=1}^{2^{d-1}} w_k \bm{a}_k
 \end{equation}
 where $\bm{a}_k = ( a_{I,k} : I \in \evenpowerset)$ for $k 
 \in \{1,\ldots,2^{d-1}\}$.
\end{theorem}
\begin{proof}
  For a vector $\bm{U} \sim C$ we can write, for any $I\in
  \mathcal{P}(\mathcal{D})$ with $|I| \geq 1 $,
\begin{align*}
     \kappa_I  = \kappa(C_I) = 2\P\left(\bm{U}_I <
    \bm{U}^*_I\right) &=  \P\left(\bm{U}_I <
    \bm{U}^*_I\right) + \P\left(\bm{U}_I >
    \bm{U}^*_I\right) \\
                        &= \P\Big(\operatorname{sign}(\bm{U}_I^* -
    \bm{U}_I)= \bm{1}\Big) + \P\Big(\operatorname{sign}(\bm{U}_I^* -
    \bm{U}_I)= -\bm{1}\Big),
\end{align*}
where $\bm{U}^*$ is an independent copy of $\bm{U}$.
If we define the random vectors
$
 \bm{V} = \operatorname{sign}(\bm{U}^* - \bm{U})$ and  $\bm{Y} = (1/2)(\bm{V} + \bm{1})$,
then the concordance probabilities of $C$ are given by
\begin{equation}\label{eq:35}
      \kappa_I  = \P(\bm{Y}_I = \bm{0}) + \P(\bm{Y}_I = \bm{1}) 
    \end{equation}
    so that the concordance signature of $C$ is determined by the
    distribution of $\bm{Y}$, which is a radially symmetric Bernoulli distribution. 
The radial symmetry follows from the fact that $\bm{V} \stackrel{d}{=} -\bm{V}$ so it must be the
case that $2\bm{Y} - \bm{1}\stackrel{d}{=}\bm{1}-2\bm{Y}$ or,
equivalently, $\bm{Y} \stackrel{d}{=} \bm{1} - \bm{Y}$.

From~\eqref{eq:35} we can conclude that, for every $I \in \powerset$,
\begin{align*}
   \kappa_I 
  &= \sum_{k=1}^{2^{d-1}} \Big( \P(\bm{Y} = \bm{s}_k) + \P(\bm{Y} =
    \bm{1} - \bm{s}_k) \Big)\indicator{I \subseteq
    J_k\,\text{or}\,I \subseteq J_k^\complement} 
= \sum_{k=1}^{2^{d-1}} w_k a_{I,k} = \kappa_I(C^*),
\end{align*}
where we have used the fact that the union of all the disjoint events
$\{\bm{Y} = \bm{s}_k\}$ and
$\{\bm{Y} = \bm{1} -\bm{s}_k\}$ forms a partition of $\{0,1\}^d$ in
the first equality, the notation~\eqref{eq:6} and~\eqref{eq:41} in the
second and~\eqref{eq:15} in
the final equality.
 Thus the weights $w_k$ specifying the extremal mixture copula are
 precisely the probabilities that specify the law of the
 radially symmetric Bernoulli vector $\bm{Y}$ through $w_k =
 2\P(\bm{Y} = \bm{s}_k)$ and the uniqueness of the set of weights
 follows from the uniqueness of the law of $\bm{Y}$.
The sufficiency of solving the
 linear equation system~\eqref{eq:34} to determine the weights $w_k$
 and the existence of a unique solution
 follow from Proposition~\ref{prop:main-result-prelim}, in particular
 the fact that~\eqref{eq:34} is a system of $2^{d-1}$ equations with
 $2^{d-1}$ unknowns and the vectors $\bm{a}_k$ are linearly independent.
\end{proof}

The implication of Theorem~\ref{theorem:conc-sign-copul} is that we can find the mixture
weights $\bm{w}$
for a given concordance signature using simple linear algebra. To see
this, let $\bm{\kappa} = \bm{\kappa}(C)$ denote the even concordance signature of a $d$-dimensional copula $C$ and let 
 $A_d$ be the $2^{d-1}\times 2^{d-1}$ matrix with columns
$\bm{a}_k$. Proposition~\ref{prop:main-result-prelim} implies that $A_d$ is of full rank, and hence invertible. The linear equation system~\eqref{eq:34} can be
written in the form $\bm{\kappa} = A_d \bm{w}$ and must have a unique
solution which can be found by
calculating $\bm{w} = A_d^{-1}\bm{\kappa}$. Theorem~\ref{theorem:conc-sign-copul} further guarantees that $\bm{w} \ge \bm{0}$ and that its components sum up to $1$. For
example, when $d = 4$ we would have
\begin{equation}\label{eq:19}
  \overbrace{\left(\begin{array}{l} 1 \\ \kappa_{\{1,2\}}
           \\\kappa_{\{1,3\}} \\\kappa_{\{1,4\}} \\
           \kappa_{\{2,3\}} \\\kappa_{\{2,4\}} \\\kappa_{\{3,4\}} \\
                     \kappa_{\{1,2,3,4\}} \end{array} \right)}^{\bm{\kappa}} =
               \overbrace{\left(
       \begin{array}{cccccccc}
         1 & 1 & 1 & 1 & 1 & 1 & 1 & 1 \\
         1 & 1 & 1 & 1 & 0 & 0 & 0 & 0 \\
         1 & 1 & 0 & 0 & 1 & 1 & 0 & 0 \\
         1 & 0 & 1 & 0 & 1 & 0 & 1 & 0 \\
         1 & 1 & 0 & 0 & 0 & 0 & 1 & 1 \\
         1 & 0 & 1 & 0 & 0 & 1 & 0 & 1 \\
         1 & 0 & 0 & 1 & 1 & 0 & 0 & 1 \\
         1 & 0 & 0 & 0 & 0 & 0 & 0 & 0
                                     \end{array} \right)}^{A_4} \bm{w}                     
                               \end{equation}

 The fact that the even concordance signature is required to determine
 the weights of the extremal mixture copula in all possible cases allows us to view the even concordance signature of a copula as the
 \textit{minimal complete} concordance signature.
 Indeed, the fact that $A_d$ is of full rank has the following
 interesting corollary, namely that a formula such as
 \eqref{eq:oddkappa} can only hold for sets of odd cardinality. 
\begin{corollary}\label{cor:kappa-formulasl}
 When $I$ is a set of even cardinality 
there is no linear formula valid for all copulas relating $\kappa_I$ to concordance
probabilities of order $|I^\prime | \leq |I|$.
                                 \end{corollary}
                                 \begin{proof}
If the statement were not true, the rows of the matrix $A_d$ would be linearly dependent,
contradicting the assertion that $A_d$ is of full rank. 
\end{proof}
\begin{remark}
 The implications of
 Corollary~\ref{cor:kappa-formulasl} can be extended.
 In view of~\eqref{eq:14} an analogous result could
 be stated for the multivariate Kendall's tau coefficients: when $|I|$
 is odd, a linear
 formula relating $\tau_I$ to the values for lower-dimensional subsets
 exists (see Proposition~1 in
\citet{bib:genest-neslehova-ben-ghorbal-11}) but no such formula
exists when $|I|$ is even. Moreover, as we will discuss in
Section~\ref{sec:main-result-example}, the concordance probabilities
of elliptical copulas are equal to twice the orthant probabilities of
Gaussian distributions centred at the origin: thus if $\bm{Z} \sim
\mathcal{N}(\bm{0}, \Sigma)$ is a
Gaussian random vector, recursive linear
formulas exist for the orthant probabilities $\P(\bm{Z}_I \geq \bm{0})$ when
$|I|$ is odd but not when $|I|$ is even.
\end{remark}


\section{Concordance signature estimation}\label{sec:signature-estimation}

Consider a random sample $\boldsymbol{X}_1,\ldots, \boldsymbol{X}_n$
from a distribution with copula $C$ and continuous margins
$F_1,\ldots, F_d$. In this section, we explain how the full concordance
signature $\tilde{\bm{\kappa}}(C)$ can be estimated intrinsically,
i.e.,~in such a way that the estimated signature is attainable.
We do this under the assumption that there are no ties in the data; this is not restrictive because the continuity of the margins ensures the absence of ties with probability~$1$.

For any $I$ with $|I| \ge 2$, empirical estimators of
$\kappa_I= \kappa(C_I)$ can be derived from empirical estimators of
$\tau_I = \tau(C_I)$ using~\eqref{eq:14}. When $d=2$ and $I = \{k,\ell\}$ for some distinct $k, \ell \in \{1,\ldots, d\}$, the classical estimator of $\tau_I$ going back to~\cite{Kendall:1938} and~\cite{Hoeffding:1947}~is
$$
 {\tau}_{I,n} = -1+ \frac{4}{n(n-1)} \sum_{ i \neq j} \indicator{X_{ik} \le X_{j k} , X_{i\ell} \le X_{j\ell}}\;.
$$
This is a special case of the estimator of $\tau_I$ for $| I| \geq 2$ proposed and investigated by \cite{bib:genest-neslehova-ben-ghorbal-11}, which is given by 
$$
{\tau}_{I,n} = \frac{1}{2^{\abs{I}-1} -1} \Bigl\{ -1 + \frac{2^{\abs{I}}}{n(n-1)} \sum_{ i \neq j} \prod_{k \in I} \indicator{X_{ik} \le X_{j k}}  \Bigr\}\;.
$$
Plugging this estimator into~\eqref{eq:14} yields an empirical
estimator of $\kappa_I$ of the form
\begin{equation*}
{\kappa}_{I,n} = \frac{2}{n(n-1)} \sum_{ i \neq j} \prod_{k \in I} \indicator{X_{ik} \le X_{j k}}.
\end{equation*}
From the theory of U-statistics \citep{Hoeffding:1948}, we know that the empirical concordance signature $\tilde{\bm{\kappa}}_{n} = (\kappa_{I,n} :  I \in \powerset)$ satisfies 
    $
    \sqrt{n}(\tilde{\bm{\kappa}}_{n}-\tilde{\bm{\kappa}})
   \rightsquigarrow \mathcal{N}(\boldsymbol{0}, \Sigma)
   $ 
   as $n \to
    \infty$, where $\tilde{\bm{\kappa}}=\tilde{\bm{\kappa}}(C)$, $\Sigma$ is the covariance matrix of the
    random vector  with components $C_I(\boldsymbol{U}_I) + \bar
    C_I(\boldsymbol{U}_I)$ and $\bar C_I$ is the survival function of
    $C_I$.
The following result  shows that the empirical concordance signature $\tilde{\bm{\kappa}}_n$ is in fact the concordance signature of a $d$-dimensional copula.  

\begin{theorem}
Assuming that $n \ge 2$ and there are no ties in the sample, there
exists a  unique $d$-dimensional extremal mixture copula $C_n$ such that $\tilde{\bm{\kappa}}_n= \tilde{\bm{\kappa}}(C_n)$.
\end{theorem}

\begin{proof} 
Let $\bm{Y}_{ij}
 = \big(\sign(\bm{X}_i - \bm{X}_j)+1\big)/2$ for $i \neq j$ and set
 \begin{equation}\label{eq:31}
   \widehat{w}_k = \frac{2}{n(n-1)} \sum_{i < j}\left(
   \indicator{\bm{Y}_{ij} = \bm{s}_k} + \indicator{\bm{Y}_{ij} =
     \bm{1} - \bm{s}_k} \right).
 \end{equation}
These empirical weights are estimators of $w_k = \P(\bm{Y} = \bm{s}_k) +\P(\bm{Y} = \bm{1} - \bm{s}_k)$, the weights of the extremal mixture which has the same concordance signature as $C$. Because the sample is assumed to have no ties, the weights $\widehat{w}_k$ are clearly positive and sum to $1$. Thus they describe a mixture of extremal copulas, say $C_n$. To establish the result, it suffices to show that for any $I\in \mathcal{P}(\mathcal{D})$ with $|I|\geq 2$, $\kappa_{I,n} = \kappa(C_{n,I})$, where $C_{n,I}$ is the margin of $C_n$ corresponding to the index set $I$. This can be seen as follows. First,
\begin{align*}
\kappa(C_{n,I})=  \sum_{k=1}^{2^{d-1}} \widehat{w}_k
a_{I,k} &=  \sum_{k=1}^{2^{d-1}}  \frac{2}{n(n-1)} \sum_{i < j}\left(
   \indicator{\bm{Y}_{ij} = \bm{s}_k} + \indicator{\bm{Y}_{ij} =
          \bm{1} - \bm{s}_k} \right) a_{I,k}\\
 & = \frac{2}{n(n-1)} \sum_{i < j}  \sum_{k=1}^{2^{d-1}} \left(
   \indicator{\bm{Y}_{ij} = \bm{s}_k} + \indicator{\bm{Y}_{ij} =
          \bm{1} - \bm{s}_k} \right) \indicator{I \subseteq
    J_k\,\text{or}\,I \subseteq J_k^\complement}.
\end{align*}    
Second, use the fact that
  \begin{align*}  
 \sum_{i < j} \left(\indicator{\bm{Y}_{ij,I} =
    \bm{0}} + \indicator{\bm{Y}_{ij,I} =
    \bm{1}} \right) &=  \sum_{i < j}\left( \indicator{\bm{X}_{i,I} <
    \bm{X}_{j,I}} + \indicator{\bm{X}_{i,I} >
    \bm{X}_{j,I}} \right) \\
    &=  \sum_{i \neq j} \indicator{\bm{X}_{i,I} <
    \bm{X}_{j,I}} = \sum_{i \neq j} \prod_{k \in I}
      \indicator{X_{ik} < X_{jk}},
\end{align*}
where the last expression is $\kappa_{I,n}  \times n(n-1)/2$ since there are no ties in the sample. 
\end{proof}

\begin{remark}
  While the probability of ties in a
  sample from a distribution with continuous margins is zero, rounding
effects may lead to occasional ties in practice. In this case it may be that some of the vectors $\bm{Y}_{ij}$ have components equal to
0.5. Let us suppose that $\bm{Y}_{ij}$ has $k$ such values for $k
\in \{1,\ldots,d\}$. A possible approach to incorporating this
information in the estimator
is to replace $\bm{Y}_{ij}$ by the $2^k$ vectors that have zeros and
one in the same positions as $\bm{Y}_{ij}$, each weighted by $2^{-k}$,
and to generalize~\eqref{eq:31} to be a weighted sum of indicators. For example, the
observation $\bm{Y}_{ij} = (1, 0, 0.5, 0.5)$ would be replaced by
$(1,0,1,1)$, $(1,0,1,0)$, $(1,0,0,1)$ and $(1,0,0,0)$, each weighted
by 1/4. This would still deliver estimates $\widehat{w}_k$ that are
positive and sum to one and thus yield a proper concordance
signature.
\end{remark}

\section{Applications}

Theorem~\ref{theorem:conc-sign-copul} allows us to test whether a vector 
of putative concordance probabilities $\bm{\kappa}_S = (\kappa_I : I
\in S)$ with label set $S$ could be a partial concordance signature of some copula $C$. We now present a number of applications of this result. Let us first say that the vector $\bm{\kappa}_S$ is 
{\it attainable} if there exists a $d$-dimensional copula $C$ such that $\bm{\kappa}_S = \bm{\kappa}_S(C)$.


\subsection{Kendall rank correlation matrices}\label{sec:Kendall}

We first turn to the characterization of matrices of pairwise Kendall's taus, also termed Kendall rank correlation matrices, which are widely used to summarize pairwise associations in a random vector. In view of~\eqref{eq:14}, this can be achieved using Theorem~\ref{theorem:conc-sign-copul} by considering the label set $S = \{\emptyset,\{1,2\},\ldots,\{1,d\},\ldots,\{d-1,d\}\}$.  

For a random vector $\bm{X}$ we denote the Kendall rank correlation matrix and the linear correlation matrix by $P_\tau(\bm{X})$ and $P_\rho(\bm{X})$ respectively.  First note that Kendall rank correlation
matrices are correlation matrices, i.e.,~positive semi-definite
matrices with ones on the diagonal and all elements in $[-1,1]$. This is because 
 if $\bm{U}$ is a random vector distributed
as a copula $C$,  then $P_\tau(\bm{U}) = 
P_\rho(\operatorname{sign}(\bm{U} - \bm{U}^*))$, where
$\bm{U}^*$ is an independent copy $\bm{U}^*$.

In Theorem~\ref{theorem:kendalls-tau1} below we provide an answer to
the attainability question for Kendall rank correlation matrices using Theorem~\ref{theorem:conc-sign-copul} and the following corollary. Recall from Section~\ref{sec:notat-gener-case} that the Kendall rank correlation  matrix, which is also the linear correlation matrix of the $k$th extremal copula, is given by $P^{(k)} = (2\bm{s}_k -\bm{1}) (2\bm{s}_k -\bm{1})^\top$.

\begin{corollary}\label{lemma:kendalls-tau}
Let $\bm{U}$ be distributed according to the mixture of extremal
copulas given by $C^* = \sum_{k=1}^{2^{d-1}} w_{k} C^{(k)}$. Then
\begin{equation}
  \label{eq:32}
  P_\tau(\bm{U}) = P_\rho(\bm{U}) = \sum_{k=1}^{2^{d-1}} w_k P^{(k)}\;.
\end{equation}
  \end{corollary}
  \begin{proof}
Proposition~\ref{prop:multi-concordance} and equation \eqref{eq:14} imply $P_\tau(\bm{U}) = \sum_{k=1}^{2^{d-1}} w_k
P_\tau(\bm{U}^{(k)})$; the equality
$ P_\rho(\bm{U}) = \sum_{k=1}^{2^{d-1}} w_k P_\rho(\bm{U}^{(k)})$ was
proved by~\cite{bib:tiit-96}. Hence~\eqref{eq:32} follows from the equality of linear and Kendall correlations for
extremal copulas.
    \end{proof}

\begin{theorem}\label{theorem:kendalls-tau1}
The  $d\times d$ correlation matrix $P$ is a Kendall rank correlation matrix
if and only if $P$ can be represented as a convex combination of the
extremal correlation matrices in dimension $d$, that is,
\begin{equation}
  \label{eq:21}
  P = \sum_{k=1}^{2^{d-1}} w_k P^{(k)}\;.
\end{equation}
\end{theorem}
\begin{proof}
If $P$ is of the form~\eqref{eq:21} then
 it is the Kendall's tau
matrix of the extremal mixture copula $C^* = \sum_{k=1}^{2^{d-1}} w_{k} C^{(k)}$ by Corollary~\ref{lemma:kendalls-tau}. If $P$ is the
Kendall's $\tau$ matrix of an arbitrary copula $C$ then, by Theorem~\ref{theorem:conc-sign-copul}, it is also the
Kendall's tau matrix of the extremal mixture copula with the same
concordance signature and
must take the form~\eqref{eq:32} by Corollary~\ref{lemma:kendalls-tau}.
\end{proof}

The set of convex combinations of extremal correlation matrices is
known as the cut polytope.
\citet{bib:laurent-poljak-95} showed that the cut polytope is a strict
subset of the so-called elliptope of correlation matrices in dimensions $d \geq
3$; see
also Section 3.3 of~\citet{bib:hofert-koike-19}. For
example, the positive-definite correlation matrix
\begin{displaymath}
 \frac{1}{12} \left(\begin{array}{rrr}
          12 & -5 & - 5 \\
          -5 & 12 & -5 \\
          -5 & - 5 & 12\\
        \end{array}
      \right)
    \end{displaymath}
    is in the elliptope but not in the cut polytope. Therefore it
    cannot be a matrix of Kendall's tau values. If $\tau = -5/12$ and
    we set $\kappa = (1 + \tau)/2$ and $\bm{\kappa} = (1,\kappa,\kappa,
    \kappa)$, there is no solution to the equation $A_3 \bm{w} =
    \bm{\kappa}$ on the 4-dimensional unit simplex and a
    representation of the form~\eqref{eq:21} is impossible.

\citet{bib:huber-maric-19} have shown that the cut polytope is also the
set of attainable linear correlation matrices for multivariate
distributions with symmetric Bernoulli margins. This
can be deduced from Theorem~\ref{theorem:conc-sign-copul} and
Proposition~\ref{prop:bernoulli} using the following lemma.
\begin{lemma}\label{lemma:tau-Bernoulli-rho}
Let $\bm{U} = U\bm{B} +(1-U)(\bm{1}-\bm{B})$, where $\bm{B}$ is a Bernoulli vector with symmetric
margins and $U$ is a uniform random variable independent
of $\bm{B}$. Then $ P_\tau\left(\bm{U}\right) = P_\rho(\bm{B})$.
  \end{lemma}
  \begin{proof}
    First note that $P_\tau \left( \bm{U}\right) = P_\rho\left(
      \bm{U}\right)$ by Proposition~\ref{prop:bernoulli} and Corollary~\ref{lemma:kendalls-tau}. By writing
\begin{align*}
\bm{U}\bm{U}^\top & =  U^2 \bm{B}\bm{B}^\top + (1-U)^2 (\bm{1} - \bm{B})(\bm{1} - \bm{B})^\top + U(1-U)\left(\bm{B}(\bm{1} - \bm{B})^\top + (\bm{1} - \bm{B})\bm{B}^\top \right)\\
 & =  (2U-1)^2 \bm{B}\bm{B}^\top + (2U-1)(1-U)\big( \bm{B}\bm{1}^\top + \bm{1}\bm{B}^\top\big) + (1-U)^2 \bm{1}\bm{1}^\top
\end{align*}
and using the fact that $\E(B_i) = 0.5$ for all $i \in \{1,\ldots,d\}$, we find that
\begin{align*}
  \cov(\bm{U}) &= 
\frac{1}{3}\Bigl( \E(\bm{B}\bm{B}^\top) - \frac{1}{2} \E(\bm{B}\bm{1}^\top+\bm{1}\bm{B}^\top) + \bm{1}\bm{1}^\top \Bigr)
-\frac{1}{4}\bm{1}\bm{1}^\top   = \frac{1}{3} \E(\bm{B}\bm{B}^\top) -
  \frac{1}{12} \bm{1}\bm{1}^\top\,.
\end{align*}
Since $\var(U_i) = 1/12$ and
$\var(B_i) =1/4$ for all $i \in \{1,\ldots,d\}$
we conclude that
$
P_\rho\left( \bm{U}\right) = 4 \E(\bm{B}\bm{B}^\top) -  \bm{1}\bm{1}^\top = P_\rho(\bm{B}).
$
This completes the proof.
\end{proof}


\begin{proposition}
The set of Kendall rank correlation matrices of copulas is identical to the set of linear
correlation matrices of Bernoulli random vectors with symmetric margins.
\end{proposition}
\begin{proof}
 If $P_\tau(\bm{U})$ is a Kendall rank correlation matrix 
 then, by Theorem~\ref{theorem:conc-sign-copul},
 it is identical to the Kendall rank correlation matrix of a random vector
 $\bm{U}^*$ with distribution given by the extremal mixture
 copula with the same concordance signature as $\bm{U}$. It
 follows from
 Proposition~\ref{prop:bernoulli} that $\bm{U}^*$ has the stochastic
 representation $\bm{U}^* \eqdis U\bm{Y} +(1-U)(\bm{1}-\bm{Y})$
 where, without loss of generality, $\bm{Y}$ has a radially symmetric
 Bernoulli distribution (with symmetric
 margins). Lemma~\ref{lemma:tau-Bernoulli-rho} gives $P_\tau(\bm{U}) =
P_\tau(\bm{U}^*) =P_\rho(\bm{Y})$.

 Conversely if $P_\rho(\bm{B})$ is the correlation matrix of a Bernoulli
 random vector $\bm{B}$ with symmetric margins (not necessarily radially
 symmetric) then, by Proposition~\ref{prop:bernoulli}, we can take an independent uniform random variable
 $U$ and construct a random vector $\bm{U}^* = U\bm{B}
 +(1-U)(\bm{1}-\bm{B})$ with an extremal mixture
 copula. Lemma~\ref{lemma:tau-Bernoulli-rho} gives $P_\rho(\bm{B}) = P_\tau(\bm{U}^*)$.
  \end{proof}

\begin{remark}
It is clear from the proof that the set of linear
correlation matrices of Bernoulli random vectors with symmetric
margins is equal to the set of linear
correlation matrices of radially symmetric Bernoulli random
vectors. This insight also appears in Theorem 1 of~\citet{bib:huber-maric-19}.
\end{remark}

\subsection{Attainability of concordance signatures}\label{sec:atta-conc-sign}

We now turn to the problem of determining whether a putative partial concordance signature is attainable, and, if it is, of calculating the set in which the remaining unmeasured concordance probabilities must lie. 

\begin{figure}[htb!]
  \centering
  \includegraphics[width=5cm,height=5cm]{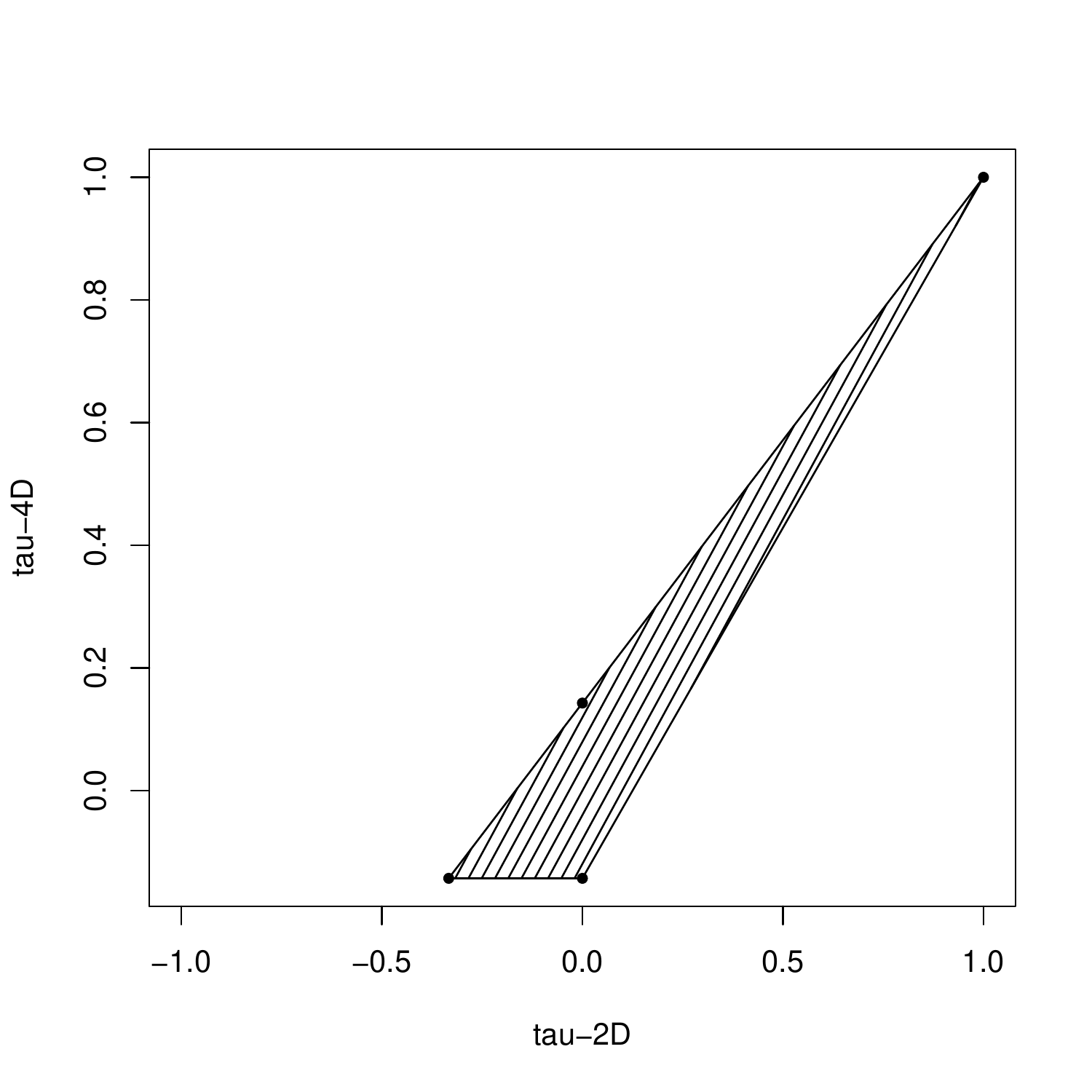}
  \includegraphics[width=5cm,height=5cm]{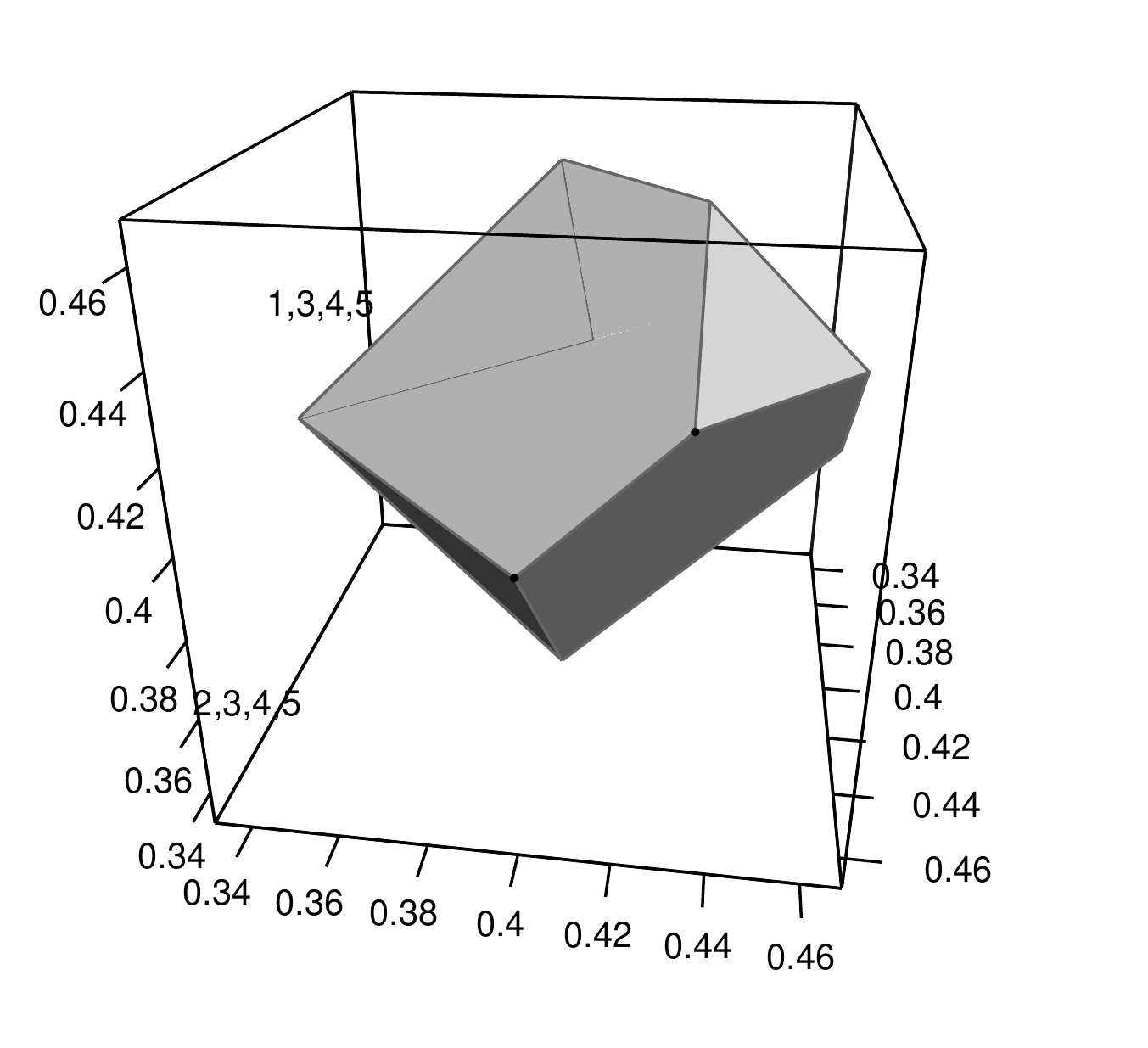}
 \caption{Left panel: set of attainable values of $(\tau_2,\tau_4)$ in Example~\ref{example1}. Right panel: Convex polytope of attainable fourth order concordance
   probabilities in Example~\ref{ex:polytope}.}
  \label{fig:polytope}
\end{figure}

We begin with a simple example that illustrates what happens in four
dimensions when all pairwise Kendall rank correlations are equal and
form a so-called equicorrelation matrix.
\begin{example}\label{example1}
  Let $C$ be a copula in dimension $d=4$ with Kendall rank correlation
  matrix $P_\tau(C)$ equal to the equicorrelation matrix with
  off-diagonal element $2\kappa_2 -1$; in other words the bivariate concordance
  probabilities $\kappa_{\{i,j\}}$ for all pairs of random variables are equal to
  $\kappa_2$. The even concordance signature is then $\bm{\kappa}(C) =
  (1,\kappa_2,\kappa_2,\kappa_2,\kappa_2,\kappa_2,\kappa_2,\kappa_4)$ where
  $\kappa_4 =
\kappa_{\{1,2,3,4\}}$.
It is easy to verify that the linear system
given in~\eqref{eq:19} is solved by the weight vector
\begin{displaymath}
  \bm{w} = \left(\kappa_4, w_1(\kappa_2,\kappa_4),
    w_1(\kappa_2,\kappa_4), w_2(\kappa_2,\kappa_4),
    w_1(\kappa_2,\kappa_4), w_2(\kappa_2,\kappa_4), w_2(\kappa_2,\kappa_4), w_1(\kappa_2,\kappa_4)\right),
\end{displaymath}
where $w_1(\kappa_2,\kappa_4) = (3\kappa_2 -1)/2-\kappa_4$ and
$w_2(\kappa_2,\kappa_4) = 1-2\kappa_2 + \kappa_4$. Since the weights
must satisfy $0 \leq w_i(\kappa_2,\kappa_4) \leq 1$, the concordance
probabilities must satisfy $\kappa_2 \in [1/3,1]$ and $\kappa_4 \in [\max(2\kappa_2 -1,0),
(3\kappa_2-1)/2]$. The left panel of Figure~\ref{fig:polytope} shows the
set of attainable values for $\tau_2 = 2\kappa_2-1$ and $\tau_4 = (8
\kappa_4 -1)/7$.  It is notable that the mere fact that the bivariate
concordance probabilities are equal limits the range of attainable
correlations significantly. In principle, for any copula, the pair $(\tau_2,\tau_4)$ 
must always lie in the rectangle $[-1,1]\times[-1/7,1]$ shown in the
plot. However, in the case of equicorrelation,
the attainable set is considerably smaller. \new{We provide some more
general results on equiconcordant copulas (copulas for which $\kappa_I
= \kappa_{\tilde{I}}$ when $|I| = |\tilde{I}|$) in Section~\ref{sec:conc-sign-equic}.}
\end{example}

In general, attainability and compatibility problems can be solved by linear
programming. To see how this is achieved, consider a label set $S \subset \evenpowerset$ strictly
contained in the even power set and suppose
we want to test whether the vector $\bm{\kappa}_S = (\kappa_I
: I \in S)$ with label set $S$ is attainable.  Let $A_d^{(1)}$ be the $|S| \times 2^{d-1}$ matrix consisting of the
  rows of $A_d$ that correspond to $S$; let $A_d^{(2)}$ be
  the matrix formed of the remaining rows of $A_d$. Consider the set
  of weight
  vectors $\{\bm{w} : A_d^{(1)} \bm{w} = \bm{\kappa}_S,\;
    \bm{w} \geq \bm{0}\}$ and note that every element of the set
    satisfies the sum condition on the weights, since
    $S$ is a label set containing $\emptyset$ and the first row of
    $A_d^{(1)}$ consists of ones. If $\bm{\kappa}_S$ is an attainable partial concordance signature
   then this set of weight vectors is non-empty and forms a convex polytope, that is, a
    set of the form
     \begin{equation}
      \label{eq:25}
        \left\{ \sum_{i=1}^m \mu_i \bm{w}_i,\;
    \sum_{i=1}^m \mu_i = 1,\; \mu_i\geq 0, i=1,\ldots,m\right\}\;.
    \end{equation}
 In this case it is possible to use the method of~\citet{bib:avis-fukuda-92} to find the \textit{vertices}
    $\bm{w}_1,\ldots,\bm{w}_m$ of the polytope. The set of attainable
    even concordance signatures containing the partial
    signature $\bm{\kappa}_S$ is then given by the convex hull of the points
    $\bm{\kappa}_i=A_d\bm{w}_i$, $i\in\{1,\ldots, m\}$, while the set
    of attainable unspecified
    elements of the concordance signature is given by the convex
    hull of the points
    \begin{equation}
      \label{eq:36}
      \{\bm{\kappa}_i^{(2)}=A_d^{(2)}\bm{w}_i,\; i = 1,\ldots,    m\}.
    \end{equation}

As the dimension increases and the discrepancy between the length of
$\bm{\kappa}_S$ and the length of an even
signature increases, the vertex enumeration algorithm can become
computationally infeasible. In such cases we may be content to simply
test $\bm{\kappa}_S$ for attainability by finding a single even concordance signature that contains $\bm{\kappa}_S$.
To do so,
 we can
  attempt to solve the optimization problem
  \begin{equation}\label{eq:23}
    \min \| A_d^{(2)} \bm{w} \| : A_d^{(1)} \bm{w} = \bm{\kappa}_S,\;
    \bm{w} \geq \bm{0}
  \end{equation}
  where $\|\cdot\|$ denotes the Euclidean norm.
  This is a standard minimization problem with both equality and
  inequality constraints. The putative partial concordance signature $\bm{\kappa}_S$ is
  attainable if a solution exists and the solution to the optimization
  problem will be the
  weight vector which gives the (collectively) smallest values for the missing
  concordance probabilities. To get the
  (collectively) largest
  values we could solve
 \begin{equation}\label{eq:24}
    \min \| A_d^{(2)} \bm{w} - \bm{1} \| : A_d^{(1)} \bm{w} = \bm{\kappa}_S,\;
    \bm{w} \geq \bm{0}.
  \end{equation}
Note also that if we are only concerned with finding the smallest or
largest missing values for certain missing concordance probabilities, then we
can drop rows of $A_d^{(2)}$.

    \begin{example}\label{ex:polytope}
For $d=5$ let a partial concordance signature be given by
$\kappa_{\{i,j\}} = 2/3$ for all pairs of variables and
$\kappa_{\{1,2,3,4\}} = \kappa_{\{1,2,3,5\}} = 0.4$. To complete the
concordance signature, three further concordance probabilities are
required: $\kappa_{\{1,2,4,5\}}$, $\kappa_{\{1,3,4,5\}}$ and
$\kappa_{\{2,3,4,5\}}$. Using the method of~\citet{bib:avis-fukuda-92}
the set of possible weight vectors is non-empty and has
nine vertices; thus the specified partial signature is attainable. The
set of attainable values for the missing values forms a polytope
in 3d which is shown in Figure~\ref{fig:polytope}.
\end{example}

\subsection{Data illustration} \label{sec:data}

We now return to the motivating example at the 
beginning of the paper. Using real data on cryptocurrency returns we illustrate the use of the
signature estimation method
and we show how missing
rank correlation values may be inferred from existing values.
\begin{example}\label{example-data}
 We take the multivariate time series of cryptocurrency prices (in US
  dollars) for Bitcoin, Ethereum, Litecoin and Ripple. From these data
  we compute the daily log-returns for the calendar year 2017, giving
  us 365 4-dimensional data points. The estimated even concordance
  signature is
  \begin{displaymath}
    \bm{\kappa}_n = (1, 0.639, 0.666, 0.598, 0.681, 0.630, 0.661, 0.364)
  \end{displaymath}
  while the weight vector describing the extremal
  mixture copula $C$ satisfying $\bm{\kappa}(C) = \bm{\kappa}_n$ is
  \begin{displaymath}
    \bm{w} = (0.364, 0.129, 0.069, 0.077, 0.098, 0.075, 0.066, 0.122),
  \end{displaymath}
  where all figures are given to 3 decimal places. The Kendall rank
  correlations of the copula $C$ will exactly match the estimated
 values from the data.

  Suppose we are not interested in Ripple but rather in a new
  cryptocurrency X4coin for which we have no data.
  We want to estimate the correlations between X4coin and the first three cryptocurrencies. By computing the convex hull of
  the set of points in~\eqref{eq:36}, and converting the concordance
  probabilities to Kendall's tau values, we obtain the polytope shown in the left
  panel of Figure~\ref{fig:estimation}. 

  Suppose that we now decide that a plausible value for
  $\tau_{\{1,4\}}$ is 0.598, which is actually the estimated rank
  correlation between Bitcoin and Ripple implied by the estimated concordance signature
above. Then the remaining two values must lie in the convex set shown
in the right panel of Figure~\ref{fig:estimation}, which is a section of the 3-dimensional
set, shown as a cut in the left panel. This is certainly the case for the
estimated values of the Etherium-Ripple and Litecoin-Ripple rank
correlations which are shown
as a point within the set.
\end{example}

\section{\new{Concordance signatures of elliptical copulas}}\label{sec:main-result-example}

The concordance signature is identical for the copulas
of all continuous elliptical
distributions with the same correlation matrix $P$. This follows
because the individual probabilities of concordance are identical for
all such copulas. This is proved in~\citet[Section
2.1]{bib:genest-neslehova-ben-ghorbal-11} in the context of an
analysis of multivariate Kendall's tau coefficients.

Let
$\bm{X}$ have an elliptical distribution centred at the origin with
dispersion matrix equal to the
correlation matrix $P$ and assume that $\P(\bm{X} = \bm{0}) = 0$. If $\bm{X}^*$ is an
independent copy of $\bm{X}$, then the concordance probabilities~\eqref{eq:1}
are given by
$\kappa_{I} = 2\P\left( \bm{X}_I < \bm{X}_I^*\right) = 2\P\left(
  \bm{X}_I - \bm{X}_I^* < \bm{0}\right)$. The random vector  $\bm{X}_I
- \bm{X}_I^*$ also has an elliptical distribution centred at the
origin. Using the stochastic representation for elliptical
distributions \citep{bib:fang-kotz-ng-90,bib:mcneil-frey-embrechts-15} we can write $\bm{X}_I
- \bm{X}_I^* \eqdis R_1 A \bm{S}$ and $\bm{X}_I \eqdis R_2 A \bm{S}$
where $\bm{S}$ is a random vector uniformly distributed on the unit sphere,
$A$ is a matrix such that $AA^\top = P$ and $R_1$ and $R_2$ are
positive scalar random variables, both independent of
$\bm{S}$. It follows that
\begin{equation}\label{eq:26}
  \kappa_{I} = 2\P\left(
  \bm{X}_I - \bm{X}_I^* < \bm{0}\right)  = 2 \P\Big((R_1 A \bm{S})_I < \bm{0}\Big) = 2\P \Big((R_2 A \bm{S})_I <
  \bm{0}\Big) = 2\P(\bm{X}_I < \bm{0}).
\end{equation}

From this calculation, we see that the positive scalar random variable
$R_2$ plays no role, so that the concordance probabilities $\kappa_I$
are the same for any elliptical random vector with the same dispersion
matrix $P=AA^\top$. Moreover, they are equal to twice the orthant
probabilities for a centred elliptical distribution with dispersion
matrix $P$. In practice it is easiest to calculate the orthant probabilities of a
multivariate normal distribution $\bm{X} \sim \mathcal{N}(\bm{0}, P)$ and this
is the approach we take in our examples; see~\ref{app:B} for an illustration involving a $6 \times 6$ correlation matrix.

Every linear correlation matrix can be the correlation matrix of a multivariate
elliptical (or multivariate normal) distribution. We now show by
means of a counterexample that an analogous statement is not true of
Kendall rank correlation matrices.

\begin{example}
The
positive-definite correlation matrix
\begin{equation}
  \label{eq:28}
 P_\tau =  \left(
    \begin{array}{rrrr}
      1 & - 0.19 & -0.29 & 0.49 \\
      -0.19 & 1 & -0.34 & 0.30 \\
      -0.29 & -0.34 & 1 & -0.79 \\
      0.49 & 0.30 & -0.79 & 1
      \end{array}
    \right)
  \end{equation}
  is a Kendall rank correlation matrix but is not the Kendall rank correlation matrix of an
  elliptical distribution.
  The elements of this matrix are attainable
  values for the Kendall's tau coefficients $\tau_{\{i,j\}}$ if  the corresponding
  concordance probabilities $\kappa_{\{i,j\}} = (1 +
  \tau_{\{i,j\}})/2$ are attainable.
  Using the methods of the previous section we can verify that
 $$
 \bm{\kappa} =
  (\kappa_\emptyset, \kappa_{\{1,2\}}, \kappa_{\{1,3\}}, \kappa_{\{1,4\}},
  \kappa_{\{2,3\}}, \kappa_{\{2,4\}}, \kappa_{\{3,4\}} )
  $$
  is an
  attainable partial concordance signature. Solving the linear
  programming problem~\eqref{eq:23} gives the weight vector
  \begin{displaymath}
    \bm{w}_1 = (0.04, 0.005,
  0.36, 0, 0.0625, 0.2475, 0.2825, 0.0025)
  \end{displaymath}
  corresponding to the minimum attainable
  fourth order concordance probability of 0.04. Solving the linear
  programming problem~\eqref{eq:24} gives the weight vector
  \begin{displaymath}
    \bm{w}_2 = (0.0425, 0.0025,
  0.3575, 0.0025, 0.06, 0.25, 0.285, 0)
  \end{displaymath}
  corresponding to the maximum attainable
  fourth order concordance probability of 0.0425.
  In this case $\bm{w}_1$ and $\bm{w}_2$ are precisely the two
  vertices of the polytope of attainable weights given
  by the set~\eqref{eq:25}, which takes the form of a line segment
  connecting $\bm{w}_1$ and $\bm{w}_2$. Any weight vector in this set will
  give the Kendall rank correlation matrix $P_\tau$.
  
 Now let us assume that $P_\tau$ in~\eqref{eq:28} corresponds to an elliptical
  copula.
  \cite{bib:lindskog-mcneil-schmock-03} and~\citet{bib:fang-fang-02} have shown that the 
 Kendall rank correlation matrix of an elliptical copula with
 correlation matrix $P$ is given by the componentwise transformation $P_\tau = 
 2\pi^{-1}\arcsin(P)$. It must be the case that $P = \sin (\pi P_\tau/2)$ is the correlation
  matrix of the elliptical copula. However, by calculating the eigenvalues we
  find that $P$ is not positive semi-definite, which is a
  contradiction.
  \end{example}

  We now turn to the copula of the multivariate Student $t$
  distribution $C^t_{\nu,P}$ with degree of freedom parameter $\nu$
  and correlation matrix parameter $P$. It is unusual to consider this copula for degrees of
  freedom $\nu <1 $, but \new{as it turns out, the Student $t$ copula converges to an extremal mixture as $\nu
  \to 0$}. 
The proof of the following remarkable property relies
 on some limiting results for the univariate and multivariate $t$
 distribution which are collected in~\ref{app:C}.

\begin{theorem}\label{thm:spider}
As $\nu\to 0$ the $d$-dimensional $t$ copula $C^t_{\nu,P}$ converges pointwise to the
unique extremal mixture copula that shares its concordance signature.
  \end{theorem}

\begin{proof}
Let the function $h_\nu(w,\bm{s})$ be defined by
\begin{equation*}
  h_\nu(w,\bm{s}) = F_\nu\left(G_{d,\nu}^{-1}(w) A \bm{s}
  \right),\quad w \in (0,1),\quad \bm{s} = (s_1,\ldots,s_d),\;
  \bm{s}^\top \bm{s} =1,
\end{equation*}
where $F_\nu$ is the distribution function of a $t$ distribution with $\nu$ degrees of
freedom, $G_{d,\nu}$ is the distribution function of the radial component of a
$d$-dimensional multivariate $t$ distribution with $\nu$ degrees of
freedom and $A$ is a $d\times d$ matrix such that $A
  A^\top =P$; such a matrix can be constructed for any positive
  semi-definite $P$.
Let $\bm{S}=(S_1,\ldots,S_d)$ be uniformly
distributed on the unit sphere and let $W$ be an independent uniform
random variable. Then $\bm{X} = G_{d,\nu}^{-1}(W)
A \bm{S}$ has a multivariate $t$ distribution and
$\bm{U} = h_\nu(W, \bm{S}) $ has joint distribution function $C_{\nu.P}^t$. We want to
show that the joint distribution function of $h_\nu(W, \bm{S})$ converges to the joint
distribution function of an extremal mixture as $\nu \to 0$.

We first argue that the random vector given by
  $A\bm{S}$ satisfies $\P\big( (A \bm{S})_j = 0\big) = 0$. Let $\bm{a}_j$ denote
the $j$th row of $A$. If $\bm{a}_j = \bm{0}$ then the $j$th row and
column of $P$ would consist of zeros implying that the $j$th margin of
the multivariate $t$ distribution of $\bm{X}$ was degenerate; this case
can be discounted because the $t$
copula with such a matrix $P$ is not defined.
Suppose therefore that $\P(\bm{a}_j^\top \bm{S} =0) >
0$ for $\bm{a}_j \neq \bm{0}$.
If $R$ is the radial random variable corresponding to the multivariate
normal, then $\P(\bm{a}_j^\top R\bm{S} =0)=\P(\bm{a}_j^\top
\bm{Z} =0) > 0$, where $\bm{Z}$ is a vector of $d$ independent standard
normal variables. However $\bm{a}_j^\top \bm{Z}$ is univariate
normal with variance $\bm{a}_j^\top \bm{a}_j > 0$ and cannot have an
atom of mass at zero.

We can therefore define the set $\mathcal{A} = \{ \bm{s} : \bm{s}^\top \bm{s} = 1, (A\bm{s})_j
\neq 0, j=1, \ldots,d\}$ such that $\P(\bm{S} \in \mathcal{A}) = 1$.
Given that $\bm{S} = \bm{s} \in \mathcal{A}$, then
\begin{align*}
  \{U_j \leq u_j \} &=
    \left\{ W \leq
  G_{d,\nu}\left(\frac{F_\nu^{-1}(u_j)}{(A\bm{s})_j}\right) \right\}\;
                      \text{if $(A\bm{s}) > 0$,} \\
   \{U_j \leq u_j \} &=
  \left\{ W \geq G_{d,\nu}\left(\frac{F_\nu^{-1}(u_j)}{(A\bm{s})_j}\right) \right\}\;
 \text{if $ (A\bm{s}) < 0$,}
\end{align*}
  and hence the conditional distribution function of $\bm{U}$ given $\bm{S}=\bm{s}$ has the form 
\begin{align*}
\P(U_1 \le u_1, \ldots, U_d \le u_d | \bm{S} = \bm{s})
& =  \P\left(\max_{j \not\in I_{A\bm{s}}} \left\{G_{d,\nu}\left(\frac{F_\nu^{-1}(u_j)}{(A\bm{s})_j}\right)\right\} \le W \le \min_{j \in I_{A\bm{s}}} \left\{G_{d,\nu}\left(\frac{F_\nu^{-1}(u_j)}{(A\bm{s})_j}\right)\right\}\right)\\
& = \left( \min_{j \in I_{A\bm{s}}} \left\{G_{d,\nu}\left(\frac{F_\nu^{-1}(u_j)}{(A\bm{s})_j}\right)\right\} -\max_{j \not\in I_{A\bm{s}}} \left\{G_{d,\nu}\left(\frac{F_\nu^{-1}(u_j)}{(A\bm{s})_j}\right)\right\} \right)^+,
\end{align*}
where $I_{A\bm{s}}$ is the set of indices $j$ for which $(A\bm{s})_j >0$. Writing, for any $\bm{u} \in (0,1)^d$,
$$
C_{\nu.P}^t(\bm{u}) = \P(U_1 \le u_1, \ldots, U_d \le u_d) = \int \P(U_1 \le u_1, \ldots, U_d \le u_d | S = \bm{s})  d F_{\bm{S}}(\bm{s}),
$$
we can use Proposition~\ref{lem:auxiliary} in~\ref{app:C} and Lebesgue's Dominated Convergence Theorem to conclude that $C_{\nu.P}^t(\bm{u})$ converges, as $\nu \to 0$, to 
\begin{equation}\label{eq:29}
 C(\bm{u}) = \int   \left( \min_{j \in I_{A\bm{s}}} (2u_j-1)^+ -\max_{j \not\in I_{A\bm{s}}} (1-2u_j)^+ \right)^+  d F_{\bm{S}}(\bm{s}).
\end{equation}
We now show that this limit is a mixture of extremal
copulas. To this end, consider  the random vector
$\big(\sign(A\bm{S})+\bm{1}\big)/2$. This has the same distribution as the multivariate
Bernoulli random vector $\bm{Y}$ whose distribution is defined by the
probabilities $p_I
= \P(\bm{Y}_I = \bm{1}) =
\P\big( (A\bm{S})_j > 0, j \in I\big)$ for $I \subseteq
\mathcal{D}$; the random vectors $\big(\sign(A\bm{S})+\bm{1}\big)/2$ and $\bm{Y}$
differ only on the null set where components of $A\bm{S}$ are zero.
Moreover, the distribution of $\bm{Y}$ is radially
symmetric since the spherical symmetry of $\bm{S}$ implies $A\bm{S}
\eqdis -A\bm{S}$ which in turn implies
$\bm{Y} \eqdis \bm{1} - \bm{Y}$.
The limiting distribution~\eqref{eq:29} may be written in the form
\begin{displaymath}
 C(\bm{u}) = \sum_{\bm{y} \in \{0,1\}^d} \left( \min_{j : y_j  = 1}
    (2u_j-1)^+ -\max_{j : y_j = 0} (1-2u_j)^+ \right)^+
  \P(\bm{Y} = \bm{y})
\end{displaymath}
and using the index set notation defined in
Section~\ref{sec:notat-gener-case} this may also be written as
\begin{align*}
 C(\bm{u}) &= \sum_{k=1}^{2^{d-1}} \left( \min_{j \in J_k}
    (2u_j-1)^+ -\max_{j \in J_k^\complement} (1-2u_j)^+ \right)^+
   \P(\bm{Y} = \bm{s}_k) \\
   &+ \sum_{k=1}^{2^{d-1}} \left( \min_{j \in J_k^\complement}
    (2u_j-1)^+ -\max_{j \in J_k} (1-2u_j)^+ \right)^+
   \P(\bm{Y} = \bm{1} - \bm{s}_k) \;.
\end{align*}
Setting $ \P(\bm{Y} = \bm{s}_k) =  \P(\bm{Y} = \bm{1}
- \bm{s}_k) = 0.5 w_k$ as in Section~\ref{sec:notat-gener-case} we obtain
$C(\bm{u}) = \sum_{k=1}^{2^{d-1}}w_k C_k(\bm{u})$, where
\begin{displaymath}
  2 C_k(\bm{u}) = \left( \min_{j \in J_k}
    (2u_j-1)^+ -\max_{j \in J_k^\complement} (1-2u_j)^+ \right)^+ + \left( \min_{j \in J_k^\complement}
    (2u_j-1)^+ -\max_{j \in J_k} (1-2u_j)^+ \right)^+
\end{displaymath}
and we need to check that $C_k$ is in fact the $k$th extremal copula $C^{(k)}$ given by~\eqref{eq:4}.
To do so, we have to distinguish the four cases below:
\begin{enumerate}
\item[(i)] Suppose that there is at least one $j \in J_k$ and at least one $j \in J_k^\complement$ such that $u_j \le 0.5$. Then $C_k(\bm{u}) = 0 = C^{(k)}(\bm{u})$.

\item[(ii)] Suppose that for all $j \in J_k$, $u_j > 0.5$ and there exists at least one $j \in J_k^\complement$ such that $u_j \le 0.5$. Then $2 C_k(\bm{u})$ equals
$$
 \left( \min_{j \in J_k}( 2u_j-1) -\max_{j \in
    J_k^\complement}(1-2u_j) \right)^+ = 2 \left( \min_{j \in J_k}u_j
  +\min_{j \in J_k^\complement}u_j -1 \right)^+ = 2C^{(k)}(\bm{u}).
$$

\item[(iii)] The case when for all $j \in J_k^\complement$, $u_j > 0.5$ and there exists at least one $j \in J_k$ such that $u_j \le 0.5$ is analogous to case (ii) and is omitted.

\item[(iv)]  Suppose that for all $j \in \mathcal{D}$, $u_j > 0.5$. Then 
$$
2 C_k(\bm{u})=\min_{j \in J_k}( 2u_j-1)  + \min_{j \in
  J_k^\complement}( 2u_j-1)  = 2\left(\min_{j \in J_k}u_j +\min_{j \in
    J_k^\complement}u_j -1\right) = 2C^{(k)}(\bm{u}).
$$
\end{enumerate}
Finally, we need to verify that the concordance signature of the limiting mixture
of extremal copulas is the same as the concordance signature of
$C^t_{\nu,P}$ for any $\nu > 0$. If $ \kappa_I$ denotes a concordance
probability for the $t$ copula, we need to show that
  $\kappa_I = \sum_{k=1}^{2^{d-1}} w_k a_{I,k}$, which is the
  corresponding concordance probability for the limit. Recall that the
  vector $\bm{X} = G_{d,\nu}^{-1}(W)
A \bm{S}$ has a multivariate $t$ distribution.
Equation~\eqref{eq:26} implies that
\begin{displaymath}
     \kappa_I = 2\P\left( G_{d,\nu}^{-1}(W) (A\bm{S})_I < \bm{0} \right)= 2\P\Big((A\bm{S})_I < \bm{0}\Big) =
    2\P(\bm{Y}_I = \bm{0}) =  2\P(\bm{Y}_I = \bm{1}) = \sum_{k=1}^{2^{d-1}} w_k a_{I,k},
  \end{displaymath}
where the final step uses the reasoning employed in
the proof of Theorem~\ref{theorem:conc-sign-copul}. 
\end{proof}

\begin{figure}[t!]
  \centering
  \includegraphics[width=8cm,height=8cm]{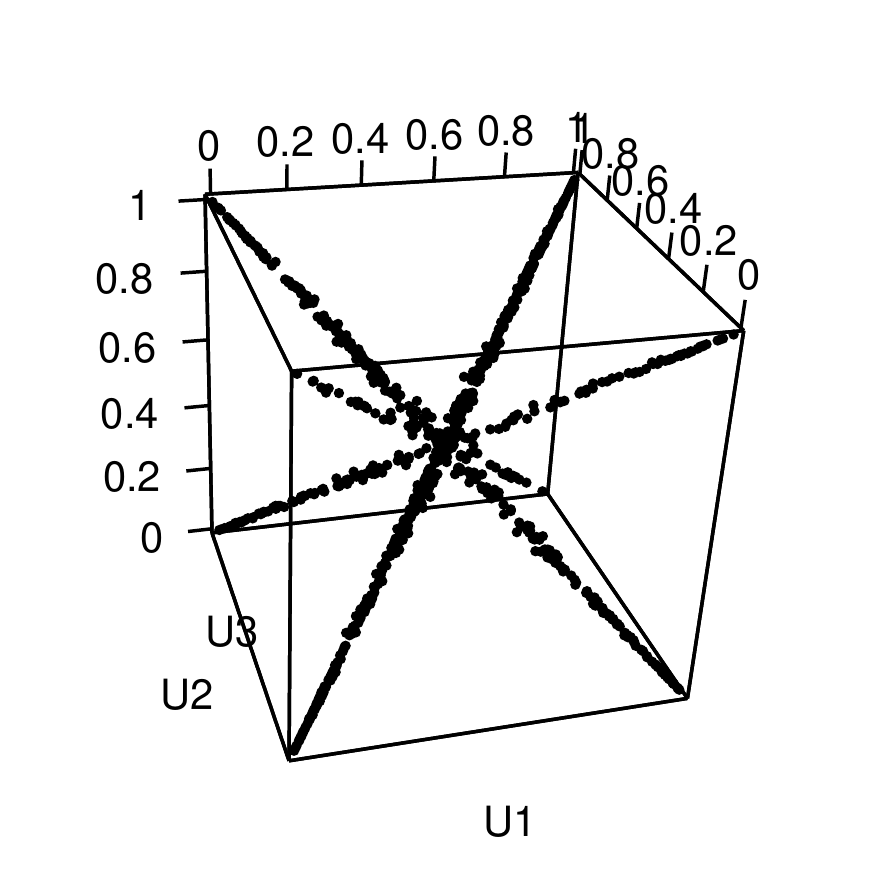}
 \caption{Scatterplot of data with distribution
$C^t_{\nu,P}$ when $\nu=0.03$ and $P$ is the $3\times3$ matrix with
elements $\rho_{12}=0.2$, $\rho_{13}=0.5$ and $\rho_{23}=0.8$.}
  \label{fig:1}
\end{figure}

Figure~\ref{fig:1} shows a scatterplot of the copula
$C^t_{\nu,P}$ when $\nu=0.03$ and $P=(\rho_{ij})$ is the $3\times3$ matrix with
elements $\rho_{12}=0.2$, $\rho_{13}=0.5$ and $\rho_{23}=0.8$. Clearly
the points are distributed very close to the four diagonals of the unit
cube. The limiting weights attached to extremal copulas associated with
each diagonal are given in Table~\ref{table:1}.
\begin{remark}
The proof holds even when the matrix $P$ is
not of full rank. However, because in such cases the copula is distributed on a
strict subspace of the unit hypercube $[0,1]^d$, the 
limiting extremal mixture copula has zero mass on certain diagonals of
the hypercube.
Suppose, for example, that rows $i$ and $j$ of the matrix $A$ satisfying $AA^\top = P$
are identical. Then the components $Y_i$ and $Y_j$ of the vector
$\bm{Y}$ defined in the proof are identical, i.e.,~they are both $0$ or
both $1$. For any vector $\bm{s}_k$ such that $s_{k,i} \neq s_{k,j}$
it must be the case that $w_k = \P(\bm{Y} = \bm{s}_k) + \P(\bm{Y} =
\bm{1} -\bm{s}_k)= 0$ and so the $k$th diagonal would have zero mass.
  \end{remark}


\begin{table}
  \caption{ \label{table:1} Probabilities associated with diagonals of
    unit cube for copula in Figure~\ref{fig:1}.}
\centering
\begin{tabular}{cr} \toprule
  Diagonal & Probability \\ \midrule
$(0,0,0) \leftrightarrow (1,1,1)$ & 51.3\%\\
$(0,0,1) \leftrightarrow (1,1,0)$ & 5.1\%\\
$(0,1,0) \leftrightarrow (1,0,1)$ & 15.4\%\\
$(0,1,1) \leftrightarrow (1,0,0)$ & 28.2\% \\ \bottomrule
\end{tabular}
\end{table}

\section{\new{Concordance signatures of equiconcordant copulas}}\label{sec:conc-sign-equic}

As we saw in Example~\ref{example1}, the constraint that certain concordance probabilities are equal further restricts the set of attainable signatures. We close this paper by investigating this phenomenon in greater detail. To this end, let $\Pi(\bm{u})$ represent a permutation of the vector $\bm{u} = (u_1,\ldots,u_d)$ and recall that a copula $C$ is said to be exchangeable if $C(\Pi(\bm{u})) = C(\bm{u})$ for all $\bm{u} \in [0,1]^d$ and any permutation $\Pi$. A weaker notion of symmetry based on concordance signatures is the  following:
\begin{definition}
A copula $C$ is said to be equiconcordant if its even concordance signature
$\bm{\kappa}(C) = \left\{ \kappa_I : I \in \evenpowerset \right\}$ has the property that $\kappa_I = \kappa_{\tilde{I}}$ whenever $|I| = |\tilde{I}|$.
\end{definition}
Every exchangeable copula $C$ is equiconcordant but not every equiconcordant copula is exchangeable. However, as we will show, the notions are equivalent for the class of extremal mixture copulas.

To develop the necessary arguments we introduce some further notation for extremal copulas building on the notation of Section~\ref{sec:notat-gener-case}. Let $\eta_k = \max(|J_k|,|J_k^\complement|)$ define the \textit{comonotonic number} of the extremal copula $C^{(k)}$ or, in other words, the size of the larger of the two groups of comonotonic random variables described by $C^{(k)}$. If we order the distinct values of $\eta_k$ from largest to smallest we obtain a vector $(h_1,\ldots,h_{m(d)}) = (d, d-1, \ldots, \lceil d/2 \rceil)$ with $m(d) = 1 + \lfloor d/2 \rfloor$ elements, where $\lfloor \cdot \rfloor$ and $\lceil \cdot \rceil$ are the floor and ceiling functions respectively.

Each value $h_i$ can be associated with a \textit{multiplicity}, which is the number of extremal copulas with comonotonic number $h_i$. Provided $d \neq 2 h_i$, the multiplicity is the number of ways of choosing the $h_i$ members of the larger comonotonic group from a set of $d$ variables and is therefore $\mu_i = \tbinom{d}{h_i}$. When $d = 2 h_i$ the binomial coefficient double counts the number of extremal copulas with comonotonic number $h_i$; for example the extremal copula for which the first $h_i$ variables are comonotonic is identical to the extremal copula for which the last $h_i$ variables are comonotonic. Thus, in this case, the mutiplicity is given by $\mu_i = \tfrac{1}{2}\tbinom{d}{h_i}$.

To understand the notation it is helpful to consider a concrete case. When $d=4$ we have $\eta_1 = 4$, $\eta_2 = \eta_3 = \eta_5 = \eta_8 = 3$ and $\eta_4 = \eta_6 = \eta_7 = 2$; thus in this case the 8 extremal copulas can be split into $m(4) = 3$ groups according to comonotonic number and we find $(h_1,h_2,h_3) = (4,3,2)$ and $(\mu_1,\mu_2,\mu_3) = (1,4,3)$.

We also need to consider the behaviour of extremal copulas under
permutations of the arguments. Clearly, the first extremal copula
$C^{(1)}$ is exchangeable, but the others are not. 
For example, under the permutation
$\Pi(\bm{u}) = (u_3,u_4,u_1,u_2)$ we find that
$$
C^{(2)}(\Pi(\bm{u})) = (\min(u_3, u_4, u_1) + u_2 -1)^+
= (\min_{j \in J_5} u_j + \min_{j \in J_5^\complement} u_j - 1)^+ = C^{(5)}(\bm{u}).
$$
and, more generally, $C^{(2)} \to C^{(5)} \to C^{(2)}$, $C^{(3)} \to
C^{(8)} \to C^{(3)}$, $C^{(4)} \to C^{(4)}$, $C^{(6)} \to C^{(6)}$ and
$C^{(7)} \to C^{(7)}$. It is clear that every extremal copula is mapped either to
itself or to another extremal copula with the same comonotonic
number. For other permutations the mapping may differ and extremal
copulas may form cycles of length greater than two, but they will
always be mapped to copulas with the same comonotonic number.



\begin{proposition}
Let $C^* = \sum_{k=1}^{2^{d-1}} w_k C^{(k)}$ be an extremal mixture copula. The following statements are equivalent:
\begin{enumerate}
    \item $C^*$ is exchangeable.
    \item $C^*$ is equiconcordant.
    \item For all $k,l \in \{1,\ldots, 2^{d-1}\}$ the weights satisfy $w_k = w_l$ whenever $\eta_k = \eta_l$.
\end{enumerate}
\end{proposition}

\begin{proof}
We show that $1 \Rightarrow 2 \Rightarrow 3 \Rightarrow 1$. The first of these implications is immediate. If $I$ and $\tilde{I}$ are two subsets of $\mathcal{D}$ with the same cardinality then exchangeability implies $C^*_I =C^*_{\tilde{I}}$  and therefore $\kappa_I = \kappa(C^*_I) = \kappa(C^*_{\tilde{I}}) = \kappa_{\tilde{I}}$.

$2 \Rightarrow 3$. We use an inductive argument based on comonotonic
number $h_i$ for $i
\in\{1,\ldots,m(d)\}$. Observe first that for $i=1$ the set $I =
\mathcal{D}$ is the single subset of $\mathcal{D}$ with cardinality
$|I| = d = h_1$ and in this case $\kappa_I = w_1$, the weight on the
first extremal copula $C^{(1)}$, which is the only extremal copula
with comonotonic number $h_1$.  Now suppose that $w_k = w_l$ whenever
$\eta_k = \eta_l = h_j$ for all $j \in \{1,\ldots,i\}$ and $i <
m(d)$. We want to show that $w_k = w_l$ whenever $\eta_k = \eta_l = h_{i+1}$.
Consider the subsets $I$ such that $|I| = h_{i+1}$. For each of these
subsets we have a formula for the concordance probability of the form
$\kappa(C^*_I) = \sum_{k=1}^{2^{d-1}} w_k a_{I,k}$ by
Proposition~\ref{prop:multi-concordance} and if $a_{I,k} = 1$
it must be the case that $\eta_k \geq h_{i+1}$ by \eqref{eq:6}. For each
distinct subset $I$
there is an equal number of copulas $C^{(k)}$ with $a_{I,k} = 1$ and
$\eta_k = h_j$ for each $j \in \{1,\ldots,i\}$ and, by assumption, the
weights on copulas with the same comonotonic number are
equal. If $h_{i+1} \neq d/2$ there are $\mu_{i+1}$ extremal
copulas $C^{(k)}$ with $\eta_k = h_{i+1}$ and $\mu_{i+1}$ sets $I$ with $|I| = h_{i+1}$ and for
each distinct set $I$ there is a distinct extremal copula $C^{(k)}$
such that $a_{I,k} = 1$. If $h_{i+1} = d/2$ there are $\mu_{i+1}$
extremal copulas with $\eta_k = h_{i+1}$ and $2\mu_{i+1}$ sets $I$
with $|I| = h_{i+1}$; in this case each distinct set $I$ is weighted
on a single extremal copula but for each $C^{(k)}$ there are two sets
$I$ with $a_{I,k} = 1$. In either case, the equality of the concordance probabilities for the sets $I$ implies that the weights $w_k$ are identical for all copulas with $\eta_k = h_{i+1}$.

$3 \Rightarrow 1$. For a permutation $\Pi$ let the function $l_\Pi(k)$
give the identity of the extremal copula to which $C^{(k)}$ is mapped
under $\Pi$ and recall that this copula has the same comonotonic
number as $C^{(k)}$. Then, for all $\bm{u} \in [0,1]^d$,
$$
C^*(\Pi(\bm{u})) = \sum_{k=1}^{2^{d-1}} w_k C^{(k)}(\Pi(\bm{u})) =
\sum_{k=1}^{2^{d-1}} w_k C^{(l_\Pi(k))}(\bm{u}) = \sum_{k=1}^{2^{d-1}} w_{l_\Pi(k)} C^{(l_\Pi(k))}(\bm{u})$$
where the final step follows because $w_k = w_l$ whenever $\eta_k =
\eta_l$. To complete the proof we need to show that the function
$l_\Pi$ is simply a permutation of the indices
$k=1,\ldots,2^{d-1}$ implying that $C^*(\Pi(\bm{u})) =
C^*(\bm{u})$. To see this assume that $l_\pi(k_1) = l_\Pi(k_2)$ for
$k_1 \neq k_2$ so that, for any $\bm{u} \in [0,1]^d$,
\begin{displaymath}
  C^{(k_1)}(\Pi(\bm{u})) = C^{l(k_1)}(\bm{u}) = C^{l(k_2)}(\bm{u}) =  C^{(k_2)}(\Pi(\bm{u})).
\end{displaymath}
If $\Pi^{-1}$ denotes the inverse permutation satisfying
$\Pi(\Pi^{-1}(\bm{u})) = \bm{u}$, then this would imply that
\begin{displaymath}
  C^{(k_1)}(\bm{u}) = C^{l(k_1)}(\Pi^{-1}(\bm{u})) =
  C^{l(k_2)}(\Pi^{-1}(\bm{u})) =  C^{(k_2)}(\bm{u})
\end{displaymath}
contradicting the assumption that $k_1 \neq k_2$.
\end{proof}

This result allows us to conclude that for any equiconcordant copula $C$, the extremal mixture copula $C^*$ sharing its signature with $C$ is exchangeable. Moreover, we can infer that the set of attainable even equiconcordance signatures takes the form
$$
\left\{ \sum_{k=1}^{2^{d-1}} w_k \bm{a}_k : w_k \geq 0, \sum_{k=1}^{2^{d-1}}w_k = 1, w_k = w_l\;\text{if}\;\eta_k = \eta_l,\forall k,l \in\{1,\ldots,2^{d-1}\} \right\}
$$
and we can write this in a simplifed way as
$$
\left\{ \sum_{i=1}^{m(d)} v_i \bm{b}_i : v_i \geq 0, \sum_{i=1}^{m(d)} v_i = 1\right\},\quad
\bm{b}_i = \frac{1}{\mu_i}\sum_{k=1}^{2^{d-1}} \bm{a}_k \Ind{\eta_k = h_i}
$$
where we note that the vectors $\bm{b}_i$ are obtained by averaging the vectors $\bm{a}_k$ corresponding to extremal copulas with the same comonotonic number and are linearly independent (because the vectors $\bm{a}_k$ are linearly independent). The weights $v_i$ now express the total weight assigned to all copulas with the same comonotonic number. The set of attainable equiconcordance signatures is clearly also a convex polytope with vertices given by the vectors $\bm{b}_1,\ldots,\bm{b}_{m(d)}$.

The $2^{d-1} \times m(d)$ matrix with columns $\bm{b}_i$ has identical rows corresponding to sets $I$ and $\tilde{I}$ with identical cardinality. Duplicate rows can simply be dropped to retain exactly $m(d)$ rows corresponding to each distinct even cardinality in the set $E = \{0,2,\ldots,2\lfloor d/2 \rfloor\}$. This results in an invertible square matrix $B_d$. Let us write the elements of the even concordance signature corresponding to each distinct cardinality as $\bm{k} = (\kappa_{i} : i \in E)$ and refer to $\bm{k}$ as the skeletal signature of an equiconcordant copula. We can then consider the linear equation system $\bm{k} = B_d \bm{v}$ to solve attainability and compatibility problems for skeletal signatures.

\begin{figure}
    \centering
    \includegraphics[width=0.55\textwidth]{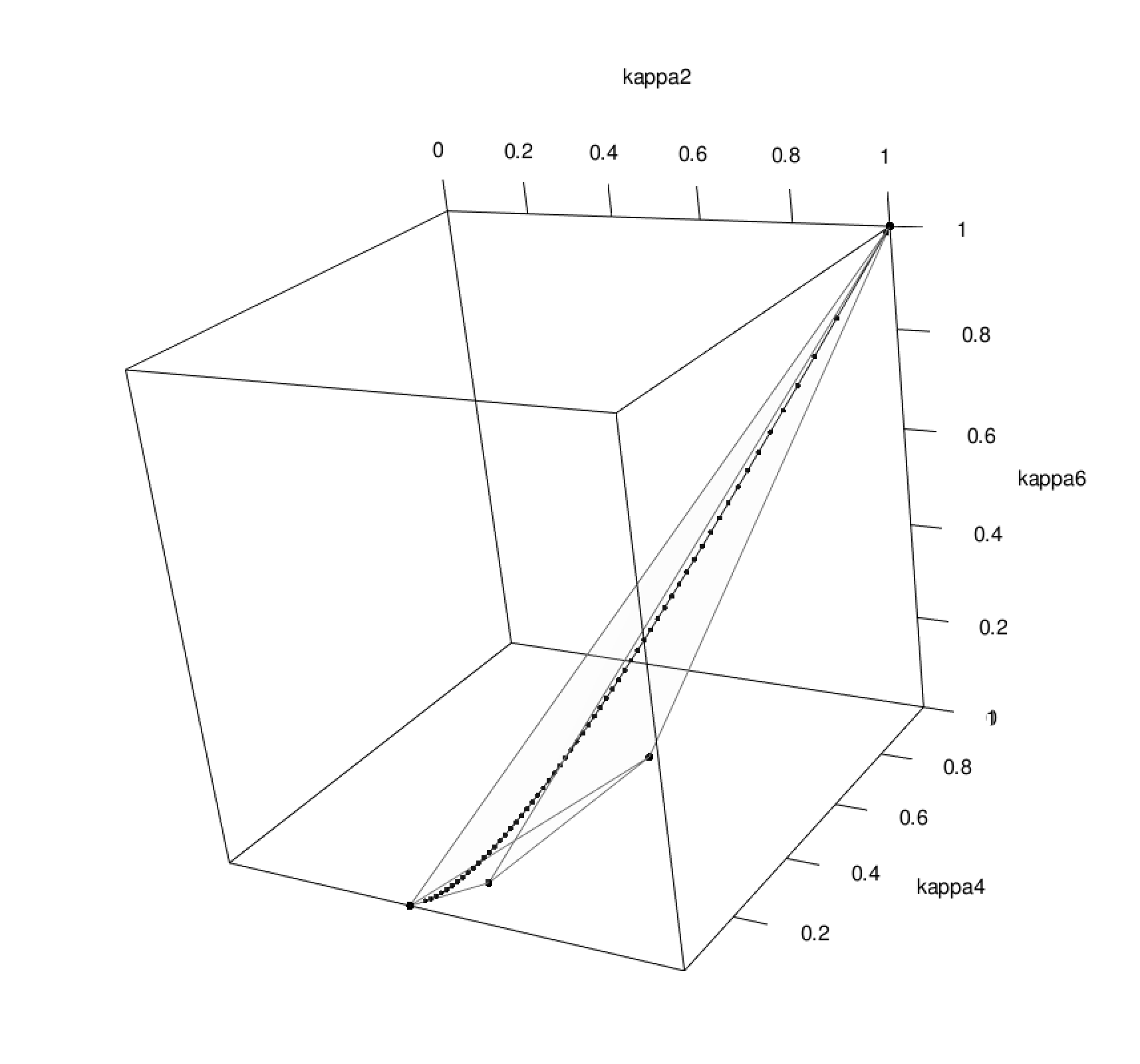}
    \caption{The tetrahedron is the set of attainable skeletal concordance signatures of equiconcordant copulas when $d=7$. The curve and associated points within the tetrahedron show the attainable skeletal signatures of exchangeable elliptical copulas.}
    \label{fig:equitope}
\end{figure}

For example, the equation system for $d=7$ is
$$
\overbrace{
\left(
\begin{array}{c}
1\\
\kappa_2\\
\kappa_4\\
\kappa_6
\end{array}
\right)}^{\bm{k}}
=
\overbrace{
\left(
\begin{array}{cccc}
1 & 1 & 1 & 1 \\[3pt]
1 & \frac{5}{7} & \frac{11}{21} & \frac{15}{35} \\[3pt]
1 & \frac{3}{7} & \frac{3}{21} & \frac{1}{35} \\[3pt]
1 & \frac{1}{7} & 0 & 0
\end{array}
\right)}^{B_7}
\overbrace{
\left(
\begin{array}{c}
v_1\\
v_2\\
v_3\\
v_4
\end{array}
\right)}^{\bm{v}}
$$
and the columns of $B_7$ form the 4 vertices of the convex polytope of attainable vectors for the skeletal signature; upon removal of the 1's in the first component these describe an irregular tetrahedron inside the unit cube $[0,1]^3$ as shown in Figure~\ref{fig:equitope}.

In higher dimensions it quickly becomes computationally infeasible to calculate the matrix $B_d$ by collapsing the $2^{d-1} \times 2^{d-1}$ matrix $A_d$. We end by giving a general formula for the elements of the matrix $B_d$ using combinatorial arguments.

\begin{theorem}
The elements of the matrix $B_d$ satisfy $B_d(1,j)=1$ for $j \in
\{1,\ldots,m(d)\}$ and
$$ B_d(i,j)=
\frac{
\left(
\begin{array}{c}
d - l_i \\
h_j - l_i
\end{array}
\right)
+
\left(
\begin{array}{c}
d - l_i \\
d - h_j - l_i
\end{array}
\right)}
{\left(
\begin{array}{c}
d  \\
h_j
\end{array}
\right)}, \quad i \in \{2,\ldots,m(d)\},\;j \in \{1,\ldots,m(d)\},
$$
where $l_i = 2(i-1)$, $h_j$ is the $j$th distinct value of the comonotonic number in reverse order from largest to smallest and, by convention, $\tbinom{n}{p} = 0$ if $p <0$.
\end{theorem}
\begin{proof}
Obviously the first row of $B_d$ consists of ones corresponding to the sum constraint on the weights so let $i > 1$.
To begin with, we exclude the case where $d$ is even and $h_j = d/2$. In all other cases the denominator of the fraction above is the multiplicity $\mu_j$ of the comonotonic number $h_j$. Let $I \subseteq \mathcal{D}$ be any subset with cardinality $|I| = l_i$.
The numerator is the number of columns of the matrix $A_d$ that correspond to extremal copulas with comonotonic number $h_j$ and which have a one in the row corresponding to $I$.

Without loss of generality let us consider the set $I$ consisting of the first $l_i$ elements of $\mathcal{D}$.
In this case we need to simply count the number of extremal copulas $C^{(k)}$ with comonotonic number $h_j$ for which the vectors $\bm{s}_k$ that we use to code the extremal copulas have zeros in the first $l_i$ entries. Now there are two possibilities: either the comonotonic number of the extremal copula $C^{(k)}$ equals the numbers of zeros in the vector $\bm{s}_k$ or the number of ones. In the first case, if there are zeros in the first $l_i$ positions of $\bm{s}_k$ there are $\tbinom{d-l_i}{h_j - l_i}$ ways of assigning the remaining zeros to the other positions. In the second case case there are $\tbinom{d-l_i}{d -h_j - l_i}$ ways of assigning the remaining zeros to the other positions. If either $h_j - l_i <0$ or $d-h_j - l_i<0$ then the corresponding binomial coefficients are simply zero.

In the case where $d$ is even and $h_j = d/2$ the denominator is twice the multiplicity $\mu_j$ of the comonotonic number $h_j$. In this case there are exactly the same number of zeros and ones in the $\mu_j$ vectors $\bm{s}_k$ coding extremal copulas with comonotonic number $h_j$ and so there are $\tbinom{d-l_i}{h_j - l_i}$ ways of assigning the remaining zeros to the other positions. However, since $\tbinom{d-l_i}{h_j - l_i}= \tbinom{d-l_i}{d -h_j - l_i}$ when $h_j =d/2$, we can use the same general formula. 
\end{proof}

\section*{Software}

The methods and examples in this paper are documented in
the \textsf{R} package \texttt{KendallSignature} at \texttt{https://github.com/ajmcneil/KendallSignature}.

\section*{Acknowledgements}
\new{The authors are grateful to the associate editor and two anonymous referees for their feedback and insightful comments that improved this paper and spurred further research that led to the material in Section~8.} This work was partially supported by a Discovery Grant from the Natural Sciences and Engineering Research Council (NSERC) of
Canada awarded to J.G. Ne\v{s}lehov\'a (RGPIN-2015-06801). 

\appendix

\section{Additional material on extremal mixture copulas}\label{app:A}

{\it Proof of  Proposition~\ref{prop:bernoulli}.}
Let $\boldsymbol{U}$ be of the form \eqref{eq:7} where $U$ and $\boldsymbol{B}$ are independent. For any $\boldsymbol{u} \in [0,1]$, 
\begin{align*}
\P(\boldsymbol{U }\le \boldsymbol{u}) & =\sum_{\boldsymbol{b} \in \{0,1\}^d} \P(\max_{j: b_j=0}(1-v_j) \le U\le \min_{j : b_j =1} v_j |\boldsymbol{B}=\boldsymbol{b})\P(\boldsymbol{B}=\boldsymbol{b})\\
& = \sum_{\boldsymbol{b} \in \{0,1\}^d} \left(\min_{j : b_j =1} v_j + \min_{j : b_j =0} v_j-1\right)^+\P(\boldsymbol{B}=\boldsymbol{b})\\
& =  \sum_{k=1}^{2^{d-1}} \left(\min_{j \in J_k} v_j + \min_{j \in J_k^\complement} v_j-1\right)^+\Big(\P(\boldsymbol{B}=\boldsymbol{s}_k)+ \P(\boldsymbol{B}=\bm{1}-\boldsymbol{s}_k)\Big).
\end{align*}
where in the final step we have used the fact that the set of possible
outcomes
of the Bernoulli vector $\bm{B}$ can be written as the disjoint union
  \begin{equation}\label{eq:18}
\{0,1\}^d =  \bigcup_{k=1}^{2^{d-1}}\{ \bm{s}_k, \bm{1} -\bm{s}_k \}.
  \end{equation}
From \eqref{eq:4}, one can see that the distribution function of
$\boldsymbol{U}$ is indeed an extremal mixture with weights as given
in~\eqref{eq:8}.

Conversely, given an extremal mixture copula of the form $C^* =
\sum_{k=1}^{2^{d-1}} w_k C^{(k)}$, it suffices to consider any Bernoulli vector
$\bm{B}$ independent of $U$ with the property that $w_k = \P(\boldsymbol{B}=\boldsymbol{s}_k)+
\P(\boldsymbol{B}=\bm{1}-\boldsymbol{s}_k)$; this is possible because the
events $\{\boldsymbol{B}=\boldsymbol{s}_k\}$ and
$\{\boldsymbol{B}=\bm{1}-\boldsymbol{s}_k\}$ are disjoint and their union
forms a partition of $\{0,1\}^d$ by~\eqref{eq:18}. We can then retrace the steps of the argument in reverse to
establish the representation~\eqref{eq:7}.
\qed

\medskip
{\it Proof of Proposition~\ref{prop:type-2-repr}.} For part (i), note that 
all probabilities of the form $p_I = \P(\bm{Y}_I = \bm{1})$ for sets $I$ with odd cardinality are fully determined by the equivalent probabilities for lower-dimensional sets of even cardinality. This follows from the fact that
\begin{displaymath}
p_I = \E\left(\prod_{i \in I}Y_i \right) = 
\E\left(\prod_{i\in I} (1 - Y_i) \right) = 1+ \sum_{A \subseteq I, |A| \ge 1} (-1)^{|A|}  \E\left(\prod_{i \in A}Y_i \right).
\end{displaymath}
When $I$ is odd, we then have 
  \begin{equation}
    \label{eq:12}
    2p_I = 1+ \sum_{A \subset I, 1 \le |A| < |I|} (-1)^{|A|} p_A = 1 - \tfrac{1}{2}|I| + \sum_{A \subset I, 2 \le |A| < |I|} (-1)^{|A|} p_A, 
  \end{equation}
  where the last equality follows from the fact that $p_A = 0.5$ whenever $|A|=1$. For example,
\begin{align*}
p_{\{1,2,3\}} &= \E\Big((1-Y_1)(1-Y_2)(1-Y_3)\Big) 
 = \frac{1}{2}\left(p_{\{1,2\}} +p_{\{1,3\}}+p_{\{2,3\}}\right) - \frac{1}{4}\;.
\end{align*}
The conclusion follows since the vector of joint event
probabilities $(p_I : I \in \powerset\setminus\emptyset)$ uniquely specifies
the distribution of a Bernoulli random vector.

For part (ii) observe that the vector $\bm{p}_{\bm{Y}}$  has length equal to
\begin{displaymath}
\sum_{j=1}^{\lfloor d/2 \rfloor} \left(\begin{array}{c} d \\
                                         2j\end{array}\right) =
                                     2^{d-1} -1.
\end{displaymath}
Radially symmetric multivariate Bernoulli distributions in
dimension $d$ have $2^{d-1} -1$ free parameters, with one
being deducted for
the sum constraint $ 2\sum_{k=1}^{2^{d-1}} \P(\bm{Y} = \bm{s}_k)
=1$. Thus the vector $\bm{p}_{\bm{Y}}$ is the minimal
vector of its kind that is required to fully specify the distribution
of $\bm{Y}$ \new{for all radially symmetric Bernoulli random vectors} $\bm{Y}$.
\qed

\begin{example}\label{ex:bivariate-ce}
  Let the joint distribution of $(U,Y_1,Y_2,Y_3)$ be specified by
  \begin{displaymath}
    \P(U \leq u, Y_1 = y_1, Y_2 = y_2, Y_3 = y_3) = \frac{1}{8}\left(
      u+ (-1)^{(y_1 + y_2 + y_3)}\frac{\theta u(1-u)}{4} \right)
  \end{displaymath}
  for $u \in [0,1]$ and $y_i \in \{0,1\}$, $i \in \{1,2,3\}$. Note
  that this is an increasing function in $u$ for any fixed $(y_1,y_2,y_3)$ and defines a valid distribution. It may be easily
  verified, by summing over the outcomes for the $Y_i$ variables, that the marginal
  distribution of $U$ is standard uniform, while letting $u \to 1$
  shows that the random vector $\bm{Y} = (Y_1,Y_2,Y_3)$ consists of iid
  Bernoulli variables with success probability $0.5$ (and is radially symmetric). It may also be
  verified that the marginal distributions $P(U \leq u, Y_i =y_i) =
  0.5u$ so that the pairs $(U,Y_i)$ are independent for all $i$, while the marginal
  distributions $P(U \leq u, Y_i =y_i, Y_j = y_j) =
  0.25 u$ so that $(U,Y_i,Y_j)$ are mutually independent for all $i
  \neq j$. 

  Now consider the vector $\bm{U} = U\bm{Y} +
  (1-U)(\bm{1}-\bm{Y})$. Since the pairs $(U,Y_i)$ are independent it
  is easy to see that the components $U_i = UY_i + (1-U)(1-Y_i)$ are
  uniform, implying that the distribution of $\bm{U}$ is a
  copula. Since the triples $(U,Y_i, Y_j)$ are mutually independent,
  the bivariate margins of $\bm{U}$ are extremal mixtures by
  Proposition~\ref{prop:bernoulli}. To calculate the copula $C$ of $\bm{U}$ we observe that, for $\bm{u}
  = (u_1,u_2,u_3)$,
  \begin{align*}
    C(u_1, u_2, u_3) &= \sum_{i=1}^4\Big(\P(\bm{U} \leq
                              \bm{u} , \bm{Y} = \bm{s}_k) + \P(\bm{U} \leq
                              \bm{u} , \bm{Y} = \bm{1} -\bm{s}_k)\Big)
    = 2\sum_{i=1}^4\P(\bm{U} \leq
                              \bm{u} , \bm{Y} = \bm{s}_k).
  \end{align*}
 The final equality follows because our model has the property that $(U,Y_1, Y_2,
  Y_3) \eqdis (1-U, 1-Y_1, 1-Y_2, 1-Y_3)$. This implies that
  $(\bm{U},\bm{Y}) = (U\bm{Y} +
  (1-U)(\bm{1}-\bm{Y}), \bm{Y}) \eqdis  (\bm{U},\bm{1} -\bm{Y})$.
 The copula $C$ has 4 distinct
 terms, each associated with a diagonal of the unit cube:
\begin{align*}
&4 C(u_1,u_2,u_3)  = \min(u_1,u_2,u_3) + \frac{\theta}{4}  \min(u_1,u_2,u_3)  \max(1-u_1,1-u_2,1-u_3)\\
& + \Ind{\min(u_1,u_2)+u_3-1\ge 0}\Big( \min(u_1,u_2) -\frac{\theta}{4} \min(u_1,u_2)\max(1-u_1,1-u_2)  -1 + u_3 + \frac{\theta}{4} u_3(1-u_3) \Big) \\
& + \Ind{\min(u_1,u_3)+u_2-1\ge 0}\Big( \min(u_1,u_3) -\frac{\theta}{4} \min(u_1,u_3)\max(1-u_1,1-u_3)  -1 + u_2 + \frac{\theta}{4} u_2(1-u_2) \Big) \\
& + \Ind{\min(u_2,u_3)+u_1-1\ge 0} \Big( \min(u_2,u_3) -\frac{\theta}{4} \min(u_2,u_3)\max(1-u_2,1-u_3) -1 + u_1 + \frac{\theta}{4} u_1(1-u_1) \Big).
\end{align*}
However, unless $\theta=0$, $C$ is not of the form \eqref{eq:4} and is not
an extremal mixture copula. 
\end{example}

{\it Proof of Proposition \ref{theorem:pairwise}.}
If $\bm{U}$ follows an extremal mixture copula then,
by Proposition~\ref{prop:bernoulli}, it has the
stochastic representation $\bm{U} \eqdis U\bm{B} + (1-U)(\bm{1} - \bm{B})$
for some Bernoulli random vector $\bm{B}$
and it is clear that all the bivariate margins have the same
structure. From the independence of $U$ and $\bm{B}$ we obtain
\begin{align*}
  \P\left(U_1 \leq u \mid \indicator{U_j = U_1}, j\neq 1\right) & =
   \P\left(U \leq u \mid \indicator{B_2 = B_1} ,\ldots, \indicator{B_d
                                                                  = B_1}\right) =\P(U\leq u) =u.
\end{align*}
Conversely, let us suppose that all the bivariate margins of $\bm{U}$
are mixtures of extremal copulas. This implies that almost surely, $\bm{U}$ takes values in the set $$
\bigcap_{i\neq j} \{ \bm{u} \in [0,1]^d :  u_j = u_i \: \text{or} \:  u_j = 1-u_i\},
$$
which simplifies to the union $\mathcal{E}
$ of the $2^{d-1}$ main diagonals of the unit hypercube, viz.
$$
\mathcal{E}= \bigcup_{k=1}^{2^{d-1}} \{\bm{u} \in [0,1]^d : u_j = u_1^{ (1-s_{k,j})}(1-u_1)^{s_{k,j}}, j\neq 1\}.
$$
Let $E_k = \{U_j = U_1^{ (1-s_{k,j})}(1-U_1)^{s_{k,j}}, j\neq 1\}$
represent the event that $\bm{U}$ lies on the $k$th diagonal and set
$w_k = \P(E_k)$. By the law of total probability, we get
$$
C(\bm{u}) = \sum_{k=1}^{2^{d-1} } w_k\P(U_1 \le u_1,\ldots, U_d \le
u_d \mid  E_k)= \sum_{k=1}^{2^{d-1} } w_k\P (U_1 \in [\max_{j \in J_k^\complement}
    (1-u_j), \min_{j \in J_k} u_j] \mid  E_k ). 
  $$
On the diagonals of the hypercube, $\{\bm{U} \in \mathcal{E}\} \cap \{ \indicator{U_j = U_1} = 1-s_{k,j}, j\neq 1\} = E_k$, so that $w_k = \P\{ \indicator{U_j = U_1} = 1-s_{k,j}, j\neq 1\}$ and conditioning on the event $E_k$ is identical to conditioning on the event $\{\indicator{U_j = U_1} = 1-s_{k,j}, j\neq 1\}$. We thus obtain that
$$
C(\bm{u}) = \sum_{k=1}^{2^{d-1} } w_k\P(U_1 \in [\max_{j \in J_k^\complement} u_j, \min_{j \in J_k} u_j] \mid  \indicator{U_j = U_1} = 1-s_{k,j}, j\neq 1) = \sum_{k=1}^{2^{d-1} } w_k C^{(k)}(\bm{u})
$$
where the last equality follows from \eqref{eq:10} and \eqref{eq:4}. This shows that $C$ is indeed a mixture of extremal copulas as claimed.
\qed

\section{Additional material on elliptical distributions}\label{app:B}

 The concordance probabilities for a continuous elliptical distribution the $6\times6$ correlation matrix
  \begin{equation}\label{eq:27}
    P = \frac{1}{16}\left(
      \begin{array}{rrrrrr}
        16 & 1 & 2 & 3 & 4 & 5 \\
        1 & 16 & 6 & 7 & 8 & 9 \\
        2 & 6 & 16 & 10 & 11 & 12 \\
        3 & 7 & 10 & 16 & 13 & 14 \\
        4 & 8 & 11 & 13 & 16 & 15 \\
        5 & 9 & 12 & 14 & 15 & 16
        \end{array}
        \right)
      \end{equation}
      are given in Table~\ref{table:2}, along with the weights of the unique mixture of extremal copulas sharing the same concordance signature.

 \begin{table}
\caption{\label{table:2} Results for the correlation matrix $P$ in~\eqref{eq:27}.}
\centering
\begin{tabular}{rlr|lr}
  \toprule
 $k$ & $S_k$ & $w_k$ & $I$ & $\kappa_I$ \\ 
  \midrule
1 & $\{0,0,0,0,0,0\}$ & 0.2627 & $\emptyset$ & 1.0000 \\ 
  2 & $\{0,0,0,0,0,1\}$ & 0.0009 & $\{1,2\}$ & 0.5199 \\ 
  3 & $\{0,0,0,0,1,0\}$ & 0.0131 & $\{1,3\}$ & 0.5399 \\ 
  4 & $\{0,0,0,0,1,1\}$ & 0.0037 & $\{1,4\}$ & 0.5600 \\ 
  5 & $\{0,0,0,1,0,0\}$ & 0.0304 & $\{1,5\}$ & 0.5804 \\ 
  6 & $\{0,0,0,1,0,1\}$ & 0.0037 & $\{1,6\}$ & 0.6012 \\ 
  7 & $\{0,0,0,1,1,0\}$ & 0.0088 & $\{2,3\}$ & 0.6224 \\ 
  8 & $\{0,0,0,1,1,1\}$ & 0.0179 & $\{2,4\}$ & 0.6441 \\ 
  9 & $\{0,0,1,0,0,0\}$ & 0.0579 & $\{2,5\}$ & 0.6667 \\ 
  10 & $\{0,0,1,0,0,1\}$ & 0.0029 & $\{2,6\}$ & 0.6902 \\ 
  11 & $\{0,0,1,0,1,0\}$ & 0.0100 & $\{3,4\}$ & 0.7149 \\ 
  12 & $\{0,0,1,0,1,1\}$ & 0.0108 & $\{3,5\}$ & 0.7413 \\ 
  13 & $\{0,0,1,1,0,0\}$ & 0.0165 & $\{3,6\}$ & 0.7699 \\ 
  14 & $\{0,0,1,1,0,1\}$ & 0.0085 & $\{4,5\}$ & 0.8019 \\ 
  15 & $\{0,0,1,1,1,0\}$ & 0.0063 & $\{4,6\}$ & 0.8391 \\ 
  16 & $\{0,0,1,1,1,1\}$ & 0.0659 & $\{5,6\}$ & 0.8869 \\ 
  17 & $\{0,1,0,0,0,0\}$ & 0.1037 & $\{1,2,3,4\}$ & 0.2804 \\ 
  18 & $\{0,1,0,0,0,1\}$ & 0.0029 & $\{1,2,3,5\}$ & 0.2977 \\ 
  19 & $\{0,1,0,0,1,0\}$ & 0.0114 & $\{1,2,3,6\}$ & 0.3150 \\ 
  20 & $\{0,1,0,0,1,1\}$ & 0.0091 & $\{1,2,4,5\}$ & 0.3244 \\ 
  21 & $\{0,1,0,1,0,0\}$ & 0.0193 & $\{1,2,4,6\}$ & 0.3437 \\ 
  22 & $\{0,1,0,1,0,1\}$ & 0.0073 & $\{1,2,5,6\}$ & 0.3675 \\ 
  23 & $\{0,1,0,1,1,0\}$ & 0.0062 & $\{1,3,4,5\}$ & 0.3702 \\ 
  24 & $\{0,1,0,1,1,1\}$ & 0.0390 & $\{1,3,4,6\}$ & 0.3909 \\ 
  25 & $\{0,1,1,0,0,0\}$ & 0.0338 & $\{1,3,5,6\}$ & 0.4161 \\ 
  26 & $\{0,1,1,0,0,1\}$ & 0.0064 & $\{1,4,5,6\}$ & 0.4581 \\ 
  27 & $\{0,1,1,0,1,0\}$ & 0.0076 & $\{2,3,4,5\}$ & 0.4503 \\ 
  28 & $\{0,1,1,0,1,1\}$ & 0.0232 & $\{2,3,4,6\}$ & 0.4725 \\ 
  29 & $\{0,1,1,1,0,0\}$ & 0.0100 & $\{2,3,5,6\}$ & 0.4993 \\ 
  30 & $\{0,1,1,1,0,1\}$ & 0.0136 & $\{2,4,5,6\}$ & 0.5427 \\ 
  31 & $\{0,1,1,1,1,0\}$ & 0.0036 & $\{3,4,5,6\}$ & 0.6153 \\ 
  32 & $\{0,1,1,1,1,1\}$ & 0.1831 & $\{1,2,3,4,5,6\}$ & 0.2627 \\ 
   \bottomrule
\end{tabular}
\end{table}

 \section{Some limiting properties of the univariate and multivariate
  Student distribution as
   $\nu \to 0$}\label{app:C}
   
\subsection{Univariate Student $t$ distribution}
 Let $F_\nu$, $F_\nu^{-1}$  and $f_\nu$ denote the df, inverse df and density of a 
 univariate $t$ distribution with $\nu$ degrees of freedom.
 The df satisfies
\begin{equation}\label{eq:30}
F_\nu(x) - 0.5= \frac{x
  \Gamma(\frac{\nu+1}{2})}{\sqrt{\pi\nu}\Gamma(\frac{\nu}{2})}
{}_2F_1\left(\frac{1}{2},\frac{\nu+1}{2}; \frac{3}{2}; -
  \frac{x^2}{\nu}\right),\quad x\in \R,\; \nu > 0,
\end{equation}
where ${}_2F_1$ denotes the hypergeometric function and $\Gamma$ the
gamma function. We will show that, for fixed $u \neq 0.5$, the quantile
function $F_\nu^{-1}(u)$ is unbounded as a function of $\nu$ as $\nu \to 0$. To that end we
first prove the following lemma.
\begin{lemma}\label{lemma:t-cdf}
$\lim_{\nu \to 0} F_\nu(x) = 0.5$ for all $x\in\R$.
\end{lemma}
\begin{proof}
  The lemma is trivially true for $x=0$ so we consider $x \neq 0$.
Making the substitution $y = x/\sqrt{\nu}$ in~\eqref{eq:30} gives
\begin{displaymath}
\frac{2F_\nu(\sqrt{\nu}y) - 1}{\nu} = \frac{2y \Gamma(\frac{\nu+1}{2})}{\sqrt{\pi}\nu\Gamma(\frac{\nu}{2})} {}_2F_1\left(\frac{1}{2},\frac{\nu+1}{2}; \frac{3}{2}; - y^2\right)
\end{displaymath}
and we can use the limits~\citep{bib:abramowitz-stegun-65}
\begin{displaymath}
\lim_{\nu \to 0} y \;{}_2F_1\left(\frac{1}{2},\frac{\nu+1}{2}; \frac{3}{2}; - y^2\right)
= y\; {}_2F_1\left(\frac{1}{2},\frac{1}{2}; \frac{3}{2}; - y^2\right) = \ln\left(y + \sqrt{1+y^2}\right)
\end{displaymath}
\begin{displaymath}
\lim_{\nu \to 0} \frac{2\Gamma(\frac{\nu+1}{2})}{\sqrt{\pi}\nu\Gamma(\frac{\nu}{2})}
=
\lim_{\nu \to 0} \frac{\Gamma(\frac{\nu+1}{2})}{\sqrt{\pi}} \frac{1}{\frac{\nu}{2}\Gamma(\frac{\nu}{2})} = \lim_{\nu \to 0}\frac{\Gamma(\frac{\nu+1}{2})}{\sqrt{\pi}} \frac{1}{\Gamma(\frac{\nu+2}{2})} = 1
\end{displaymath}
to conclude that
\begin{equation*}
  \lim_{\nu \to 0}   \frac{2F_\nu(\sqrt{\nu}y) - 1}{\nu} = \ln\left(y
    + \sqrt{1+y^2}\right) = \operatorname{sign}(y)\ln\left(|y| + \sqrt{1+y^2}\right).
\end{equation*}
Reversing the earlier substitution and setting $x = \sqrt{\nu} y$ now
gives
\begin{equation*}
 \lim_{\nu \to 0} \frac{ 2F_\nu(x) - 1}{ \nu \operatorname{sign}(x) \left( \ln \big( |x| +
    \sqrt{\nu+x^2} \big) - 0.5 \ln \nu \right)} = 1.
\end{equation*}
The result follows from the fact that the denominator tends to 0 as
$\nu \to 0$.
\end{proof}
\begin{lemma}\label{cor:t-quantile}
  \begin{displaymath}
    \lim_{\nu \to 0} F_\nu^{-1}(u) =
    \begin{cases}
      -\infty &\text{if $u < 0.5$,} \\
      0 &\text{if $u=0.5$,}\\
      \infty &\text{if $u > 0.5$.}
      \end{cases}
  \end{displaymath}
\end{lemma}
\begin{proof}
The case $u =0$ is obvious, since $F_\nu^{-1}(0.5) = 0$ for all
$\nu>0$.
To show that  $\lim_{\nu \to 0} F_\nu^{-1}(u) = -\infty$ for $u < 0.5$,
we need to show that, for all $k < 0$, there exists $\delta$ such that
$\nu < \delta \Rightarrow F_\nu^{-1}(u) \leq k$. Suppose we fix an
arbitrary $k < 0$. Since $F_\nu(k) \to 0.5$ as $\nu \to 0$, there
exists $\delta > 0$ such that $\nu < \delta$ implies that $F_\nu(k) >
u$, since $u < 0.5$. But then, for any $\nu < \delta$, it follows that
$k \in \{x : F_\nu(x) \geq u\}$ and hence $F_\nu^{-1}(u) = \inf \{x :
F_\nu(x) \geq u\} \leq k$.

Analogously, to show that  $\lim_{\nu \to 0} F_\nu^{-1}(u) = \infty$ for $u > 0.5$,
we need to show that, for all $k > 0$, there exists $\delta$ such that
$\nu < \delta \Rightarrow F_\nu^{-1}(u) \geq k$. Suppose we fix an
arbitrary $k > 0$. Since $F_\nu(k) \to 0.5$ as $\nu \to 0$, there
exists $\delta > 0$ such that $\nu < \delta$ implies that $F_\nu(k) <
u$, since $u > 0.5$. But then, for any $\nu < \delta$, it follows that
$k \notin \{x : F_\nu(x) \geq u\}$. Since $F_\nu(x)$ is a
non-decreasing function of $x$, any $y < k$ also satisfies $y \notin
\{x : F_\nu(x) \geq u\}$. Hence $F_\nu^{-1}(u) = \inf \{x :
F_\nu(x) \geq u\} \geq k$.



  \end{proof}
\subsection{Multivariate Student $t$ distribution}
If the random vector $\bm{X}$ has a 
$d$-dimensional multivariate $t$ distribution with $\nu$ degrees of
freedom, then it has the stochastic representation $\bm{X} \eqdis \bm{\mu} +
R A \bm{S}$
where $\bm{S}$ is a uniform random vector on the $d$-dimensional unit
sphere, $R$ is an independent, positive, scalar random variable such
that $R^2/d \sim F(d,\nu)$ (a Fisher--Snedecor $F$ distribution),
$\bm{\mu}$ is a location vector and $A$ is a matrix; see Section 6.3 in~\citet{bib:mcneil-frey-embrechts-15}.
Let $G_{d,\nu}$ denote the df of the radial random variable $R$ and $g_{d,\nu}$ the
corresponding density; the latter
 is given by
\begin{equation}
  \label{eq:17}
  g_{d,\nu}(r) = \frac{2}{r}\left(1 + \frac{\nu}{r^2}\right)^{-\frac{d}{2}}
  \left(1 + \frac{r^2}{\nu}\right)^{-\frac{\nu}{2}}\frac{\Gamma\left(
      \frac{\nu + d}{2}\right)}{\Gamma\left( \frac{d}{2}\right)\Gamma\left( \frac{\nu}{2}\right)}.
\end{equation}

The following limiting result is key to our analysis of the limiting
behaviour of the $t$ copula as $\nu \to 0$.  Note that the limiting
function on the right-hand side is either the df of a random variable that is uniformly distributed
 on $[0.5,1]$ or the survival function of a random variable that is
 uniformly distributed on $[0,0.5]$, depending on the
 sign of $\lambda$.
\begin{proposition}\label{lem:auxiliary}
For a constant $\lambda \neq 0$ and any $u \in [0,1]$, 
\begin{equation}
  \label{eq:A3}
  \lim_{\nu \to 0} G_{d,\nu}\left(\frac{F_\nu^{-1}(u)}{\lambda}\right)=
\begin{cases}
  (2u -1)^+& \text{if $\lambda > 0$,} \\
  (1 - 2u)^+ & \text{if $\lambda < 0$.}
  \end{cases}
\end{equation}
 \end{proposition}
\begin{proof} When $\lambda > 0$, $F_{d,\nu}(u) =
  G_{d,\nu}\left(\lambda^{-1}F_\nu^{-1}(u)\right)$ is a distribution
  function supported on $[0.5,1]$. Similarly, when $\lambda < 0$,
  $F_{d,\nu}(u) = 1-G_{d,\nu}\left(\lambda^{-1}F_\nu^{-1}(u)\right)$
  is a distribution function supported on $[0,0.5]$. 
The density of $F_{d,\nu}$ is given by
\begin{align*}
   f_{d,\nu}(u) &=
                                        g_{d,\nu}\left(\frac{|F_\nu^{-1}(u)|}{|\lambda|}\right)
                                        \frac{1}{|\lambda|}\frac{1}{f_\nu\big(F_\nu^{-1}(u)\big)}
\end{align*}
where $u \in [0.5,1)$ when $\lambda > 0$ and $u \in (0, 0.5]$
when $\lambda < 0$.

Note that, when
  $u=0.5$, $F_\nu^{-1}(u) =0$ for all $\nu$ and~\eqref{eq:A3} clearly
  holds; we thus restrict our analysis of the density to the case where $u \neq 0.5$. Using the notation
  $x_{\nu,u} = F_\nu^{-1}(u)$ and the expression~\eqref{eq:17} we
  have that
\begin{align*}
    f_{d,\nu}(u) &=  g_{d,\nu}\left(\frac{|x_{\nu,u}|}{|\lambda|}\right)
                   \frac{1}{|\lambda|}\frac{1}{f_\nu\big(x_{\nu,u}\big)}\\
  &= \frac{2}{|x_{\nu,u}|}\left(1 + \frac{\lambda^2\nu}{x_{\nu,u}^2}\right)^{-\frac{d}{2}}
  \left(1 + \frac{x_{\nu,u}^2}{\lambda^2\nu}\right)^{-\frac{\nu}{2}}\frac{\Gamma\left(
      \frac{\nu + d}{2}\right)}{\Gamma\left(
    \frac{d}{2}\right)\Gamma\left( \frac{\nu}{2}\right)}
    \frac{\sqrt{\nu\pi}\Gamma\left( \frac{\nu}{2}\right)}{\Gamma\left(
      \frac{\nu + 1}{2}\right)} \left(1 +
    \frac{x_{\nu,u}^2}{\nu}\right)^{\frac{\nu+1}{2}} \\
                 &= 2\, \frac{\sqrt{\nu + x_{\nu,u}^2}}{|x_{\nu,u}|}
                  \left(1 +
    \frac{\lambda^2\nu}{x_{\nu,u}^2}\right)^{-\frac{d}{2}}
    \left(1 + \frac{x_{\nu,u}^2}{\lambda^2\nu}\right)^{-\frac{\nu}{2}}\left(1 +
                   \frac{x_{\nu,u}^2}{\nu}\right)^{\frac{\nu}{2}}
                   \frac{\Gamma\left(
      \frac{\nu + d}{2}\right)}{\Gamma\left(
                   \frac{d}{2}\right)}
                   \frac{\sqrt{\pi}}{\Gamma\left(
                   \frac{\nu + 1}{2}\right)} \\
                 &= 2\, \underbrace{\frac{\sqrt{\nu + x_{\nu,u}^2}}{|x_{\nu,u}|}}_{1}
                   \underbrace{\left(1 +
    \frac{\lambda^2\nu}{x_{\nu,u}^2}\right)^{-\frac{d}{2}}}_{2}
  \underbrace{ \left(\frac{\lambda^2\left(\nu + x_{\nu,u}^2\right)}{\lambda^2\nu +
      x_{\nu,u}^2}\right)^{\frac{\nu}{2}}}_{3}
                   \underbrace{\frac{\Gamma\left(
      \frac{\nu + d}{2}\right)}{\Gamma\left(
                   \frac{d}{2}\right)}}_{4}
                   \underbrace{\frac{\sqrt{\pi}}{\Gamma\left(
                   \frac{\nu + 1}{2}\right)}}_{5} .
\end{align*}
The limit as $\nu \to 0$ of each of the five terms
above is 1. For terms 1, 2, 4 and 5, this is obvious from the properties of the gamma function and the fact that
$|x_{\nu,u}| \to \infty$ for $u\neq 0.5$
(Lemma~\ref{cor:t-quantile} above). For term 3 note that the term
within the outer brackets converges to $\lambda^2$ and hence the
assertion follows. The density in the limit thus satisfies
\begin{equation*}
 \lim_{\nu \to 0} f_{d,\nu}(u) =
  \begin{cases}
    2 \times\indicator{0.5 \leq u < 1} & \lambda > 0 \\
       2 \times\indicator{0 < u \leq 0.5} &  \lambda < 0.
    \end{cases}
  \end{equation*}
From Scheff\'e's Theorem we conclude that the limiting distribution of $F_{d,\nu}$ is uniform on
either $[0,0.5]$ or $[0.5,1]$, depending on the sign of $\lambda$. Hence if $\lambda > 0$, $F_{d,\nu}(u) \to (2u-1)^+$ as $\nu \to 0$, and
if $\lambda < 0$, $F_{d,\nu}(u) \to 1 - (1-2u)^+$ as $\nu \to 0$. In either
case~\eqref{eq:A3} holds.
\end{proof}



\bibliographystyle{elsarticle-harv}
\newcommand{\noopsort}[1]{}

\end{document}